\documentclass[11pt,fleqn, a4]{article}

\usepackage{inputenc} 
\usepackage[T1]{fontenc}

\usepackage{color}
\usepackage{amssymb}
\usepackage{amsmath}

\usepackage{dsfont}

\newcommand{\dcb}{\begin{array}{lll}}
\newcommand{\dce}{\end{array}}
\newcommand{\ebe}{\begin{enumerate}\setlength{\baselineskip}{13pt}\setlength{\parskip}{0pt}}
\newcommand{\dbe}{\end{enumerate}}

\newcommand{\ibegin}{\begin{itemize}\setlength{\baselineskip}{19pt}\setlength{\parskip}{7pt}}
\newcommand{\iend}{\end{itemize}}
\newcommand{\ok}{\rule{4pt}{6pt}}%

\newtheorem{Theorem}{Theorem}[section]
\newtheorem {Cor}[Theorem]{Corollary}
\newtheorem {definition}[Theorem]{Definition}
\newtheorem {Lemma}[Theorem]{Lemma}
\newtheorem {rem}[Theorem]{Remark}
\newtheorem {assumption}[Theorem]{Assumption}

\newcommand {\bd}{\begin{definition}}
\newcommand {\ed}{\end{definition}}
\newcommand {\bl}{\begin{Lemma}}
\newcommand {\el}{\end{Lemma}}
\newcommand {\bcor}{\begin{Cor}}
\newcommand {\ecor}{\end{Cor}}
\newcommand {\brem }{\begin{rem} \rm }
\newcommand {\erem }{\end{rem}}
\newcommand{\bethe}{\begin{Theorem}}
\newcommand{\ethe}{\end{Theorem}}
\newcommand {\bassumption}{\begin{assumption}}
\newcommand {\eassumption}{\end{assumption}}

\def \ind{1\!\!1\!}

\def\cro#1{\langle #1\rangle}

\newcommand{\pv}{{\Gamma\!}}

\begin{document}

\begin{center}
\textbf{\Large Optional splitting formula in a progressively enlarged filtration}

\

Shiqi Song

{\footnotesize Laboratoire Analyse et Probabilités\\
Université d'Evry Val D'Essonne, France\\
shiqi.song@univ-evry.fr}
\end{center}

\

\begin{quote}

\begin{center}\textbf{Abstract}\end{center}

Let $\mathbb{F}$ be a filtration and $\tau$ be a random time. Let $\mathbb{G}$ be the progressive enlargement of $\mathbb{F}$ with $\tau$. We study the following formula, called the optional splitting formula\!\!: For any $\mathbb{G}$-optional process $Y$, there exists an $\mathbb{F}$-optional process $Y'$ and a function $Y''$ defined on $[0,\infty]\times(\mathbb{R}_+\times\Omega)$ being $\mathcal{B}[0,\infty]\otimes\mathcal{O}(\mathbb{F})$ measurable, such that
$$
Y=Y'\ind_{[0,\tau)}+Y''(\tau)\ind_{[\tau,\infty)}.
$$
(This formula can also be formulated for multiple random times $\tau_1,\ldots,\tau_k$.) 

We are interested in this formula because of its fundamental role in many recent papers on credit risk modeling, and also because of the fact that its validity is limited in scope and this limitation is not sufficiently underlined. In this paper we will determine the circumstances in which the optional splitting formula is valid. We will then develop practical sufficient conditions for that validity. Incidentally, our results reveal a close relationship between the optional splitting formula and several measurability questions encountered in credit risk modeling. That relationship allows us to provide simple answers to these questions.

\end{quote} 

\

\textbf{Keywords\!\!:} optional process, progressive enlargement of filtration, credit risk modeling, conditional density hypothesis\\

\textbf{MSC Classification\!\!:} 60G07, 60G44, 91G40, 97M30.\\

\

\

\section{Introduction}\label{intro}

The progressive enlargement of filtration is a basic technique in the credit risk modeling. Let us recall its definition (cf. \cite{J}). Let $(\Omega,\mathcal{A},\mathbb{Q})$ be a probability space equipped with a filtration $\mathbb{F}=(\mathcal{F}_t)_{t\geq 0}$ of sub-$\sigma$-algebras in $\mathcal{A}$. We assume that $\mathbb{F}$ is right-continuous and that $\mathcal{F}_0$ contains $\mathcal{N}^{\mathcal{F}_\infty}$, where, for a $\sigma$-algebra $\mathcal{T}$ contained in $\mathcal{A}$, $\mathcal{N}^{\mathcal{T}}$ denotes the $\sigma$-algebra generated by $\{A\subset \Omega\!\!: \exists B\in\mathcal{T}, A\subset B, \mathbb{Q}[B]=0\}$. Let $\tau$ be a random time (i.e. a random variable taking values in $[0,\infty]$) on $(\Omega,\mathcal{A})$. The progressive enlargement of the filtration $\mathbb{F}$ with the random time $\tau$ is the filtration $\mathbb{G}=(\mathcal{G}_t)_{t\geq 0}$ where $$
\mathcal{G}_t=\mathcal{N}^{\sigma(\tau)\vee\mathcal{F}_\infty}\vee(\cap_{s>t}(\mathcal{F}_s\vee\sigma(\tau\wedge s))),\ t\geq 0.
$$ 
According to \cite{VZ}, $$
\mathcal{N}^{\sigma(\tau)\vee\mathcal{F}_\infty}\vee(\cap_{s>t}(\mathcal{F}_s\vee\sigma(\tau\wedge s)))=\cap_{s>t}(\mathcal{N}^{\sigma(\tau)\vee\mathcal{F}_\infty}\vee(\mathcal{F}_s\vee\sigma(\tau\wedge s))).
$$ 
So, $\mathbb{G}$ is a right-continuous filtration. We denote by $\mathcal{O}(\mathbb{G})$ (resp. $\mathcal{P}(\mathbb{G})$) the $\mathbb{G}$-optional (resp. $\mathbb{G}$-predictable) $\sigma$-algebra. We define $\mathcal{O}(\mathbb{F})$ and $\mathcal{P}(\mathbb{F})$ in a similar way. See \cite{DM, Yan}.

\subsection{Background}\label{bgd}

In this paper we are interested in the measurability relations that a class of random maps can have with respect to a $\sigma$-algebra. The issue of measurability relations has been considered fundamental from the very beginning of the theory of progressive enlargement of filtration. In \cite{barlow} where a honest time $\tau$ (see Section \ref{examples} for the definition) is considered, it is shown that the $\mathbb{G}$-progressively measurable processes can be written in term of the random intervals $[0,\tau),[\tau,\infty)$ and of $\mathbb{F}$-progressively measurable processes. As for the $\mathbb{G}$-predictable processes, they satisfy a stronger formula (cf. \cite[Proposition(5.3)]{J})\!\!: For any $Y\in\mathcal{P}(\mathbb{G})$, there exist $Y',Y''\in\mathcal{P}(\mathbb{F})$ such that
\begin{equation}
Y\ind_{(0,\infty)}=Y'\ind_{(0,\tau]}+Y''\ind_{(\tau,\infty)}.
\end{equation}
Based on these relationships, it is proved that every $\mathbb{F}$-martingale is a $\mathbb{G}$-semimartingale in the case of a honest time $\tau$.  

If $\tau$ is not a honest time, we have a less precise formula (\cite[Lemme (4.4)]{J})\!\!: for any $Y\in\mathcal{P}(\mathbb{G})$, there exist a $\mathbb{F}$-predictable process $Y'$ and a function $Y''$ defined on $[0,\infty]\times(\mathbb{R}_+\times\Omega)$ being $\mathcal{B}[0,\infty]\otimes\mathcal{P}(\mathbb{F})$ measurable, such that
\begin{equation}\label{pre-split}
Y=Y'\ind_{[0,\tau]}+Y''(\tau)\ind_{(\tau,\infty)}.
\end{equation}
In particular, $[0,\tau]\cap\mathcal{P}(\mathbb{G})=[0,\tau]\cap\mathcal{P}(\mathbb{F})$. This formula is used in various computations in the filtration $\mathbb{G}$ which vary from the predictable dual projections to the orthogonal decomposition of the family of $\mathbb{G}$-martingales stopped at $\tau$. 

More recently, in \cite{Kusuoka}, the martingale representation property in $\mathbb{G}$ is studied for a Brownian filtration $\mathbb{F}$ and a random time $\tau$ satisfying the two conditions\!\!: 
\ebe
\item[(i)]
any $\mathbb{F}$-martingale is a $\mathbb{G}$-martingale (called hypothesis$(H)$), and 
\item[(ii)] 
the $\sigma$-algebras $\mathcal{G}^\circ_t=\sigma(\tau\wedge t)\vee\mathcal{F}_t, t\geq 0,$ completed by the null sets, form a right-continuous filtration. 
\dbe
The condition (ii) is a measurability condition and it is not trivial. In general $\{\tau=t\}\notin\mathcal{G}^\circ_t$, but always $\{\tau=t\}\in\mathcal{G}^\circ_{t+}$. This question will be further examined in Section \ref{examples}.

The paper \cite{BSW} considers another filtration $\mathcal{G}^\star_t=\sigma(\{\tau\leq s\}\!\!: 0\leq s\leq t)\vee\mathcal{F}_t, t\geq 0$. The filtration $\mathbb{F}$ is supposed to be a complete Brownian filtration and the random time $\tau$ to be a Cox time, i.e.$$
\tau=\inf\{t\geq 0\!\!:\Gamma_t\geq \Xi\},
$$
where $\Gamma$ is a $\mathbb{F}$-adapted càdlàg increasing process and $\Xi$ is a strictly positive random variable independent of $\mathcal{F}_\infty$. Then, it is proved that $(\mathcal{G}^\star_t)_{t\geq 0}$ is a right-continuous filtration, and consequently $\mathcal{G}_t=\mathcal{G}^\star_t$. This result is a typical example of the problem studied in \cite{VZ}\!\!: in what circumstances does the following formula hold\!\!:$$
\mathcal{T}'\vee(\cap_{n=1}^\infty\mathcal{T}_n)
=\cap_{n=1}^\infty(\mathcal{T}'\vee\mathcal{T}_n),
$$ 
where $\mathcal{T}'$ is a $\sigma$-algebra and $(\mathcal{T}_n)_{n\geq 1}$ is an inverse filtration. This interchangeability problem of \cite{VZ} is in general a very delicate issue. See \cite{CY, ES1, ES2, handel, YY} for more information. See also Section \ref{first-results} below. 

This result of \cite{BSW} is a particular case of the following question\!\!: how can the $\sigma$-algebra $\mathcal{G}_T$, where $T$ is a $\mathbb{F}$-stopping time, be factorized in terms of $\sigma(\{\tau\leq s\wedge T\!\!:s\geq 0\})$ and of $\mathbb{F}_T$. Many works on $\mathbb{G}$ depend on that decomposition, especially when the monotone class theorem is applied on $\mathbb{G}_T$. For example, we have the identity $\mathcal{G}_\infty=\sigma(\tau)\vee\mathcal{F}_\infty$ (completed by the null sets). This is required in the paper \cite{Kusuoka} in order to obtain results on the martingale representation property in $\mathbb{G}$ under the hypothesis$(H)$. When the results in \cite{Kusuoka} are extended in \cite{JS3}, one has to work with a general $\mathbb{F}$-stopping time $T$ other than $\infty$. But usually the $\sigma$-algebra $\mathcal{G}_T$ is strictly greater than $\sigma(\{\tau\leq s\wedge T\!\!:s\geq 0\})\vee\mathcal{F}_T$. A laborious computation was necessary in \cite{JS3} to get around the gap between them. To better appreciate this idea, it is to be compared with the general equality $\mathcal{G}_{T-}= \sigma(\tau\wedge T)\vee\mathcal{F}_{T-}$ (completed with null sets), a consequence of formula (\ref{pre-split}) and of the identity $\{T\leq \tau\}=\{T= \tau\wedge T\}$.

In other respects, the work \cite{biagini} requires the following fact\!\!: for the complete natural filtration $\mathbb{F}$ of a Brownian motion $W$, which is postulated to remain a $\mathbb{G}$ martingale, for any $\mathbb{G}$-martingale $X$, there exists a $\mathbb{F}$-predictable process $J$ such that $X_\tau=J_\tau$ on $\{\tau<\infty\}$. This is equivalent to say $\{\tau<\infty\}\cap\mathcal{G}_{\tau}=\{\tau<\infty\}\cap\mathcal{G}_{\tau-}$. In general these two $\sigma$-algebras are different. The gap between such $\sigma$-algebras is the subject of several papers \cite{BEKSY, BPY, J,DE,SS}. We will come back to this question later.

\subsection{The subject of the paper}\label{subject}

Recently an optional version of formula (\ref{pre-split}) has been revealed to be fundamental in credit risk modeling with progressive enlargement of filtration\!\!: for any $\mathbb{G}$-optional process $Y$, there exist a $\mathbb{F}$-optional process $Y'$ and a function $Y''$ defined on $[0,\infty]\times(\mathbb{R}_+\times\Omega)$ being $\mathcal{B}[0,\infty]\otimes\mathcal{O}(\mathbb{F})$ measurable, such that
\begin{equation}\label{opt-split}
Y=Y'\ind_{[0,\tau)}+Y''(\tau)\ind_{[\tau,\infty)}.
\end{equation}
This formula (\ref{opt-split}) is directly or indirectly involved in numerous works (cf. \cite{BSW, biagini, BJR, CJZ, ElKJJ, jiao, KL, Kusuoka, Pham, wu}). That said, this widespread use of the formula suggests caution. Indeed, unlike formula (\ref{pre-split}), formula (\ref{opt-split}) is in general not valid. We recall the well-known example of \cite{barlow}\!\!: let $\mathbb{F}$ be the natural filtration of a Brownian motion $B$ with $B_0=0$. Let $T=\inf\{t\geq 0\!\!: |B_t|=1\}$ and $\tau=\sup\{s\leq T\!\!: B_s=0\}$. Then, $X=\ind_{[\tau,\infty)}\mbox{sign}(B_T)$ is a $\mathbb{G}$-martingale, which does not satisfy formula (\ref{opt-split}). See also \cite{DM,DM2} and \cite[Proposition (5.6)]{J} for a complementary analysis.

The aim of the present paper is to make a detailed analysis of formula (\ref{opt-split}) (as well as its extension to multiple random times), to determine the circumstances of the validity of the formula, and to find sufficient conditions for that validity.

\subsection{The plan}

We will investigate the problem under the name of optional splitting formula (abbreviated as $O\!\!S\!\!F$). Until now, for the sake of clarity, we have only mentioned the case of the filtration $\mathbb{G}$ generated by a single random time $\tau$. Actually the problem can also be formulated for multiple random times $\tau_1,\ldots,\tau_k$.

In Section \ref{first-results} we begin the investigation with a single random time $\tau$. We formally introduce the notion of optional splitting formula at $\tau$ (which is simply formula (\ref{opt-split})). We draw the first consequences of this notion. We prove in Lemma \ref{predictableSPLT} that the $\mathbb{G}$-predictable processes satisfy the $O\!\!S\!\!F$ at $\tau$, and in Theorem \ref{right-continuity} that the $O\!\!S\!\!F$ at $\tau$ entails an equality between $\mathcal{G}_t$ ($t\geq 0$) and $\sigma(\{\tau\leq s\}\!\!:0\leq s\leq t)\vee\mathcal{F}_t$ (completed by the null sets). In Corollary \ref{no-hold} we formally prove that the $O\!\!S\!\!F$ can not hold in general. 

We now set the stage for the proof of the first main result. We notice that the $O\!\!S\!\!F$ problem can not be treated by itself. It is a particular case of a broader problem. We consider the family $\mathcal{L}^o$ of $\mathbb{G}$-optional subsets $A\subset \mathbb{R}_+\times\Omega$ such that, for any $\mathbb{G}$-optional process $Y$, there exists a $\mathbb{F}$-optional process $Y'$ and a function $Y''$ defined on $[0,\infty]\times(\mathbb{R}_+\times\Omega)$ being $\mathcal{B}[0,\infty]\otimes\mathcal{O}(\mathbb{F})$ measurable, such that
\begin{equation}\label{local opt-split}
Y\ind_A=(Y'\ind_{[0,\tau)}+Y''(\tau)\ind_{[\tau,\infty)})\ind_A.
\end{equation}
We say then that the optional splitting formula at $\tau$ holds on $A$. Formula (\ref{opt-split}) is the particular case of formula (\ref{local opt-split}) when $A=\mathbb{R}_+\times\Omega$. To make the difference, we call formula (\ref{opt-split}) the global optional splitting formula. The question now becomes whether $\mathbb{R}_+\times\Omega\in\mathcal{L}^o$, or more generally, exactly which elements are contained in the family $\mathcal{L}^o$. We note that, no matter whether formula (\ref{opt-split}) holds, the family $\mathcal{L}^o$ always gives good indications of what the filtration $\mathbb{G}$ looks like.

In Section \ref{local-split} we examine the properties of $\mathcal{L}^o$. Clearly, if $A\in\mathcal{L}^o$, for any $\mathbb{G}$-optional set $B\subset A$, $B\in\mathcal{L}^o$. Also (cf. Lemma \ref{predictunion}), for any sequence $(A_i)_{i\geq 1}$ of predictable elements in $\mathcal{L}^o$, $\cup_{i\geq 1}A_i$ is again an element in $\mathcal{L}^o$. From subsection \ref{left-closed} to \ref{right-closed}, we establish conditions under which a random interval $(S,T]$, where $S,T$ are $\mathbb{G}$-stopping times, belongs to $\mathcal{L}^o$. The idea behind this consideration is that the intervals $(S,T]$ are $\mathbb{G}$-predictable sets. If some of them are in $\mathcal{L}^o$, their union is an element in $\mathcal{L}^o$, which can be vast enough to give an answer to the $O\!\!S\!\!F$ problem. The results of this section are essential for the next section.

Section \ref{suff-cond} is devoted to our first main result Theorem \ref{mrt_after_default}, which gives a sufficient condition for the $O\!\!S\!\!F$. We begin with the optional splitting formula on the random intervals $[0,\tau)$, $(\tau,\infty)$ and $[\tau,\infty)$. We show in Theorem \ref{beforedefault} that $[0,\tau)\in\mathcal{L}^o$ without any supplementary condition. The cases of $(\tau,\infty)$ and $[\tau,\infty)$ are not so easy. Inspired by \cite{SongThesis, song-local-solution} we introduce a covering condition. Then, with the results of section \ref{local-split}, we prove in Theorem \ref{mrt_after_default} that, if the covering condition holds on $(\tau,\infty)$, then $(\tau,\infty)\in\mathcal{L}^o$. If the covering condition holds on $[\tau,\infty)$, then the global optional splitting formula holds.

Despite its unusual definition, the $s\!\mathcal{H}$-measure condition is satisfied in most of examples we know in the literature and Theorem \ref{mrt_after_default} is applicable there. To illustrate this fact, in Section \ref{examples} we explain how the various works mentioned in Subsection \ref{bgd} are linked with the $O\!\!S\!\!F$ and how their results can be explained as consequences of Section \ref{suff-cond}. The key point is the so-called hypothesis$(H)$. We prove in Theorem \ref{HyH} that the $O\!\!S\!\!F$ holds, whenever the hypothesis$(H)$ is satisfied under a probability measure equivalent to $\mathbb{Q}$. We also present the example of the $\natural$-model in \cite{JS2} where the hypothesis$(H)$ has not been verified, but the $O\!\!S\!\!F$ holds. We recall the only explicit examples of \cite{barlow, J} where the $O\!\!S\!\!F$ does not hold.

In Section \ref{multi-time} we tackle the problem in its general form with multiple random times $\tau_1,\ldots,\tau_k$. It is to note that, once the case of a single random time is well understood, the case of multiple random times can naturally be dealt with by induction. The true challenge lies elsewhere. In fact, the multiplicity of random times may cause an inflation of notations in an induction argument. In Section \ref{multi-time}, we adopt a definition of the multi-time optional splitting formula which is specially formulated to adapt to the induction argument. (Of course, that definition remains equivalent to the one used in the literature (cf. \cite{Pham}).) We prove the induction procedure in Theorem \ref{GmO}. We recall the widely used density hypothesis. Thereafter, we prove our second main result Theorem \ref{densityH-splitting} which states that the multi-time $O\!\!S\!\!F$ holds, whenever the density hypothesis is satisfied. These results on $O\!\!S\!\!F$ are again extended to the case of multiple random times with marks in Section \ref{multi-time-marked} (cf. Theorem \ref{densityH-splitting-marked}). 

The present paper is motivated by the use (direct or indirect) of the $O\!\!S\!\!F$ in the papers \cite{CJZ, ElKJJ, jiao, KL, Pham}, etc.. The results in Section \ref{multi-time-marked} justify this use, due to the density hypothesis. This concludes the paper.

\section{Optional splitting formula at a random time $\tau$ with respect to $\mathbb{F}$}\label{first-results}

\subsection{Definition}

When a multivariate function is viewed as a process, the time-randomness pair $(t,\omega)\in\mathbb{R}_+\times\Omega$ is privileged. Other variables will be considered as parameters. More formally, let $E$ be a space and $Y(\theta,t,\omega)$ be a map defined on
$(\theta,t,\omega)\in E\times(\mathbb{R}_+\times\Omega)$. If $E$ is considered as a space of parameters, for $\theta\in E$, we denote by $Y(\theta)$ (resp. by $Y_t(\theta)$ for $t\in\mathbb{R}_+$) the map $(s,\omega)\rightarrow Y(\theta,s,\omega)$ (resp. $\omega\rightarrow Y(\theta,t,\omega)$). For a map $\Upsilon$ defined on $\Omega$ into $E$, $Y(\Upsilon)$ denotes the map $(s,\omega)\rightarrow Y(\Upsilon(\omega),s,\omega)$.

\bd\label{df_splitting}
We say that a $\mathbb{G}$-optional process $Y$ satisfies the optional splitting formula at $\tau$ with respect to $\mathbb{F}$, if there exists a process $Y'\in \mathcal{O}(\mathbb{F})$ and a function $Y''$ defined on $[0,\infty]\times(\mathbb{R}_+\times\Omega)$ being $\mathcal{B}[0,\infty]\otimes\mathcal{O}(\mathbb{F})$-measurable, such that$$
Y=Y'\ind_{[0,\tau)}+Y''(\tau)\ind_{[\tau,\infty)}.
$$
We will denote $\mathfrak{p}^{[0,\tau)}Y=Y'$ and $\mathfrak{p}^{[\tau,\infty)}Y=Y''$. 

We say that the $\mathbb{G}$-optional splitting formula holds at $\tau$ with respect to $\mathbb{F}$, if the above property is satisfied by any $\mathbb{G}$-optional process $Y$.

\ed

\brem\label{negl}
Let $\mathcal{N}$ denote $\mathcal{N}^{\sigma(\tau)\vee\mathcal{F}_\infty}$ (cf. Section \ref{intro} for the definition). We note that the identity in Definition \ref{df_splitting} is an indistinguishable equality with respect to $\mathcal{N}$, i.e. there exists a $\mathbb{Q}$-negligible set $A$ in $\mathcal{N}^{\sigma(\tau)\vee\mathcal{F}_\infty}$ such that the map $Y$ is identical to the map $Y'\ind_{[0,\tau)}+Y''(\tau)\ind_{[\tau,\infty)}$ on $A^c$. Note also that the component $Y''$ in Definition \ref{df_splitting} is uniquely defined only on the set $[\tau,\infty)$. The maps $\mathfrak{p}^{[\tau,\infty)}Y$ designates one such component $Y''$. This absence of uniqueness does not affect the exactness of the subsequent computations, because $\mathfrak{p}^{[\tau,\infty)}Y$ will be applied on $[\tau,\infty)$. Similar observations can be made on $\mathfrak{p}^{[0,\tau)}Y$.
\erem

\brem
The term "splitting" is twofold. It obviously means that the formula is split at the random time $\tau$. But, more importantly, it implies that the measurability of $Y''(\tau)$ is factorized into two components $\sigma(\tau)$ and $\mathcal{O}(\mathbb{F})$ ($Y''\in\mathcal{B}[0,\infty]\otimes\mathcal{O}(\mathbb{F})$). Theorem \ref{F-tau} and \ref{right-continuity} below compared to \cite{J,VZ} show that these splitting properties constitute a very strong condition on the filtration $\mathbb{G}$.
\erem

In the rest of this paper we often omit the qualifying expression "with respect to $\mathbb{F}$".

\subsection{Some consequences on $\mathcal{G}_\tau$}\label{some-consequences}

In this paper, if a map $\xi$ is measurable with respect to a $\sigma$-algebra $\mathcal{T}$, we will write $\xi\in \mathcal{T}$ and say that "$\xi$ is in $\mathcal{T}$". For any random time $R$ on $\Omega$, we denote (cf. \cite{J}) $$
\dcb
\mathcal{F}_R = \sigma\{X_R\ind_{\{R<\infty\}}\!\!:\mbox{ $X$ an $\mathbb{F}$-optional process}\},\\
\mathcal{F}_{R+} = \sigma\{X_R\ind_{\{R<\infty\}}\!\!:\mbox{ $X$ an $\mathbb{F}$-progressively measurable process}\}.\\
\dce
$$

\bl\label{F}
Let $F\in\mathcal{B}[0,\infty]\otimes\mathcal{O}(\mathbb{F})$. Then, the map $\omega\in\Omega\rightarrow
F_{\tau(\omega)}(\tau(\omega),\omega)$ is in $\mathcal{F}_\tau
$ and, for $t\geq0$, the map $\omega\in\Omega\rightarrow
F_{t}(\tau(\omega),\omega)$ is in $\sigma(\tau)\vee\mathcal{F}_t$.
\el

\noindent\textbf{Proof} By the monotone class theorem, we need only to see that the stated measurability is true for $F$ in the particular form $F_{t}(s,\omega)=h(s)f_t(\omega)$ where $h\in\mathcal{B}[0,\infty]$ and $f\in\mathcal{O}(\mathbb{F})$.\ok

\bethe\label{F-tau}
Assume the optional splitting formula at $\tau$. We necessarily have $\mathcal{F}_{\tau}=\mathcal{F}_{\tau+}=\mathcal{G}_{\tau}$.
\ethe

\noindent \textbf{Proof} We know that $\mathcal{G}_\tau$ is generated by $Y_\tau$ (cf. \cite[Corollary 3.23]{Yan}), where $Y$ runs through the family $\mathcal{O}(\mathbb{G})$. By the assumption of the optional splitting formula at $\tau$, there exists a $Y''\in\mathcal{B}[0,\infty]\otimes\mathcal{O}(\mathbb{F})$, such that $$
Y_{\tau(\omega)}(\omega) = Y''_{\tau(\omega)}(\tau(\omega),\omega).
$$
According to Lemma \ref{F}, $Y_\tau\in\mathcal{F}_\tau$, and consequently, $\mathcal{G}_{\tau}\subset \mathcal{F}_{\tau}$. The theorem is proved, because always $\mathcal{F}_{\tau}\subset \mathcal{F}_{\tau+}\subset \mathcal{G}_{\tau}$. \ok

Recall the result in \cite[Proposition (5.6)]{J}. Let $M$ be a continuous uniformly integrable $\mathbb{F}$-martingale such that $M_0=0, M_\infty\neq 0$. Let $\tau=\sup\{t\geq 0\!\!: M_t=0\}$. Then, $\mathcal{F}_\tau\neq \mathcal{F}_{\tau+}$. As a consequence, we have the following corollary\!\!:

\bcor\label{no-hold}
The optional splitting formula at $\tau$ can not hold in general.
\ecor

\subsection{Trace computation of $\sigma$-algebras}

Let $D$ be a subset of $\Omega$ and $\mathcal{T}$ be a $\sigma$-algebra on $\Omega$. We denote by $D\cap\mathcal{T}$ the family of all subsets $D\cap A$ with $A$ running through $\mathcal{T}$. If $D$ itself is an element in $\mathcal{T}$, $D\cap\mathcal{T}$ coincides with $\{A\in\mathcal{T}\!\!: A\subset D\}$. We use the symbol "+" to present the union of two disjoint subsets. For two disjoint sets $D_1,D_2$ in $\Omega$, for two families $\mathcal{T}_1,\mathcal{T}_2$ of sets in $\Omega$, we denote by $D_1\cap\mathcal{T}_1+D_2\cap\mathcal{T}_2$ the family of sets $D_1\cap B_1+D_2\cap B_2$ where $B_1\in\mathcal{T}_1,B_2\in\mathcal{T}_2$.

\bl\label{union-AT}
Let $\mathcal{T}$ and $\mathcal{T}'$ be two $\sigma$-algebras. Let $D$ be a set.
\ebe
\item[(a)]
For any set $D'$, we have $$
D\cap\mathcal{T}\subset D'\cap D\cap\mathcal{T}+D'^c\cap D\cap\mathcal{T}.
$$
If $D'\in\mathcal{T}$, we have $$
D\cap\mathcal{T} = D'\cap D\cap\mathcal{T}+D'^c\cap D\cap\mathcal{T}.
$$
\item[(b)]
Let $(A_i)_{i\geq 1}$ be a sequence of sets. Suppose that $D\cap A_i\cap\mathcal{T}\subset D\cap A_i\cap\mathcal{T}'$ for all $i\geq 1$. If $A_i\in\mathcal{T}'$ for all $i\geq 1$, we also have $$
D\cap (\cup_{i\geq 1}A_i)\cap\mathcal{T}\subset D\cap (\cup_{i\geq 1}A_i)\cap\mathcal{T}'.
$$
\dbe
\el

\textbf{Proof.} We only prove the second part of the lemma. For any $C\in\mathcal{T}$, for any $i\geq 1$, there exists $C_i\in\mathcal{T}'$ such that $D\cap A_i\cap C=D\cap A_i\cap C_i$. Let$$
B_1=A_1,\ B_k=A_k\setminus (\cup_{i=1}^{k-1}A_i).
$$ 
Then, $B_k\in\mathcal{T}'$ and $D\cap C\cap B_k=D\cap C_k\cap B_k$. Therefore,$$
\dcb
D\cap (\cup_{i\geq 1}A_i)\cap C
&=&(\cup_{i\geq 1}A_i)\cap D\cap(\cup_{i\geq 1}B_i)\cap C\\
&=&(\cup_{i\geq 1}A_i)\cap D\cap (\cup_{k\geq 1}(B_k\cap C_k)).\\
\dce
$$
This proves the result. \ok

\subsection{A strong right-continuity}

In this and the following subsections, an (in)equality between two measurable functions (resp. two random processes) is to be understood as an almost sure relation (resp. an indistinguishable relation) with respect to $\mathcal{N}$ (cf. Remark \ref{negl}). 

For two elements $a,b$ in $[0,\infty]$ we denote $$
a\nmid b=
\left\{
\dcb
a&&\mbox{ if $a\leq b$}\\
\\
\infty&&\mbox{ if $a> b$}.\\
\dce
\right.
$$

\bethe\label{right-continuity}
If the optional splitting formula holds at $\tau$, then for any $t\geq 0$, $\mathcal{G}_t=\mathcal{N}\vee\sigma(\tau\nmid t)\vee\mathcal{F}_t$.
\ethe

\textbf{Proof} Let $0\leq t<\infty$. The $\sigma$-algebra $\mathcal{G}_t$ is generated by $Y_t$ for $Y\in\mathcal{O}(\mathbb{G})$. We write the optional splitting formula$$
Y=Y'\ind_{[0,\tau)}+Y''(\tau)\ind_{[\tau,\infty)},
$$
where $Y'=\mathfrak{p}^{[0,\tau)}Y$ and $Y''=\mathfrak{p}^{[\tau,\infty)}Y$. Since $Y'_t\in\mathcal{F}_t, Y''_t(\tau)\in\sigma(\tau)\vee\mathcal{F}_t$ (Lemma \ref{F}), we have $$
\dcb
Y_t&=&Y'_t\ind_{\{t<\tau\}}+Y''_t(\tau)\ind_{\{\tau\leq t\}}\\
&\in&\{t<\tau\}\cap\mathcal{F}_t+\{\tau\leq t\}\cap(\sigma(\tau)\vee\mathcal{F}_t)\\
&=&\{t<\tau\}\cap(\sigma(\tau\nmid t)\vee \mathcal{F}_t)+\{\tau\leq t\}\cap(\sigma(\tau\nmid t)\vee\mathcal{F}_{t})\\
&=&\sigma(\tau\nmid t)\vee \mathcal{F}_t,
\dce
$$
where the last equality comes from the fact that $\{t<\tau\},\{\tau\leq t\}\in\sigma(\tau\nmid t)\vee \mathcal{F}_t$. This being true for any $Y\in\mathcal{O}(\mathbb{G})$, we conclude that $\mathcal{G}_t\subset\mathcal{N}\vee\sigma(\tau\nmid t)\vee\mathcal{F}_t$. It is actually an equality, because the inverse inclusion is always true. \ok

\brem
As a matter of fact, in the above theorem we can not replace the term $\tau\nmid t$ with $\tau\wedge t$. In general,$$
\{t<\tau\}\cap(\sigma(\tau\wedge t)\vee \mathcal{F}_t)+\{\tau\leq t\}\cap(\sigma(\tau\wedge t)\vee\mathcal{F}_{t})
\neq \sigma(\tau\wedge t)\vee \mathcal{F}_t,
$$
because $\{\tau\leq t\}$ (or more precisely $\{\tau= t\}$) is not necessarily in $\sigma(\tau\wedge t)\vee \mathcal{F}_t$. This is a potential pitfall. See \cite[Chapitre IV, n$^\circ$104]{DM} which comments on \cite{DE}. (The problem no longer arises if $\{\tau=t\}$ is negligible and if $\mathbb{F}$ is complete.) See also \cite[Chapter VI.3]{protter}.

\erem

\brem
As a consequence of Theorem \ref{right-continuity}, the filtration $(\mathcal{N}\vee\sigma(\tau\nmid t)\vee\mathcal{F}_t\!\!: t\geq 0)$ is right-continuous. According to \cite{VZ}, this right-continuity is a fairly strong condition on the pair $\tau$ and $\mathbb{F}$.
\erem

\subsection{Predictable processes}

\bl\label{predictableSPLT}
The $\mathbb{G}$-predictable processes satisfy the optional splitting formula at $\tau$.
\el

\textbf{Proof} The $\mathbb{G}$-predictable processes are generated (in the sense of the monotone class) by the processes of the form$$
g(\tau\wedge a)\ind_A\ind_{]a,b]} + \ind_B\ind_{[0]},
$$
where $g$ is a bounded Borel function, $0\leq a<b$ are real numbers, $A\in\mathcal{F}_a$ and $B\in\mathcal{G}_0$. We directly verify that the process $g(\tau\wedge a)\ind_A\ind_{]a,b]}$ satisfies the optional splitting formula. As for $\ind_B\ind_{[0]}$, by \cite[Lemme(4.4)]{J}, there exist $B',B''\in\mathcal{F}_0$ such that$$
\ind_B\ind_{[0]}
=\ind_{B'}\ind_{\{0<\tau\}}\ind_{[0]}+\ind_{B''}\ind_{\{\tau=0\}}\ind_{[0]}
=\ind_{B'}\ind_{[0]}\ind_{[0,\tau)}+\ind_{B''}\ind_{[0]}\ind_{[\tau,\infty)},
$$
i.e., $\ind_B\ind_{[0]}$ also satisfies the optional splitting formula.
The lemma can now be proved by the monotone class theorem. \ok

\

\section{Optional splitting formula on $\mathbb{G}$-optional sets}\label{local-split}

\subsection{Definition and basic properties}

\bd\label{df_splitting_local}
Let $A$ be a $\mathbb{G}$-optional set. We say that a $\mathbb{G}$-optional process $Y$ satisfies the optional splitting formula on $A$ (at $\tau$ with respect to $\mathbb{F}$), if there exists a $Y'\in \mathcal{O}(\mathbb{F})$ and a function $Y''$ defined on $[0,\infty]\times(\mathbb{R}_+\times\Omega)$ being $\mathcal{B}[0,\infty]\otimes\mathcal{O}(\mathbb{F})$ measurable, such that$$
Y\ind_{A}=(Y'\ind_{[0,\tau)}+Y''(\tau)\ind_{[\tau,\infty)})\ind_{A}
$$
(an indistinguishable relation). We will denote $\mathfrak{p}^{[0,\tau)}_AY=Y'$ and $\mathfrak{p}^{[\tau,\infty)}_AY=Y''$. 

We say that the ($\mathbb{G}$-)optional splitting formula holds on $A$ (at $\tau$ with respect to $\mathbb{F}$), if the above property is satisfied for any $\mathbb{G}$-optional process $Y$.

We denote by $\mathcal{L}^o$ the family of $A\in\mathcal{O}(\mathbb{G})$ on which the ($\mathbb{G}$-)optional splitting formula (at $\tau$ with respect to $\mathbb{F}$) holds.
\ed

Obviously, this definition coincides with Definition \ref{df_splitting} when $A=\mathbb{R}_+\times\Omega$. We will call the property in Definition \ref{df_splitting} the global optional splitting formula. Comments similar to those concerning $\mathfrak{p}^{[0,\tau)}Y$ and $\mathfrak{p}^{[\tau,\infty)}Y$ in Definition \ref{df_splitting}, can be made about $\mathfrak{p}^{[0,\tau)}_AY$ and $\mathfrak{p}^{[\tau,\infty)}_AY$. 

The following properties are direct consequences of Definition \ref{df_splitting_local}.

\bl\label{stable}
Let $A$ be a $\mathbb{G}$-optional set. Let $\mathcal{S}_A$ be the family of all $\mathbb{G}$-optional processes which satisfy the optional splitting formula on $A$. Then, $\mathcal{S}_A$ is a linear space, closed by pointwise limit, by $\inf,\max$, by product operations. 
\el

\bl\label{on-subset}
Let $A,B$ be two $\mathbb{G}$-optional sets such that $B\subset A$. Then, $A\in\mathcal{L}^o$ implies $B\in\mathcal{L}^o$.  
\el

\subsection{Optional splitting formula on predictable sets}

\bl\label{predictLSST}
Let $A$ be a $\mathbb{G}$-predictable set. Then, $A\in\mathcal{L}^o$ if and only if, for any $\mathbb{G}$-optional process $Y$, $Y\ind_A$ satisfies the global optional splitting formula.
\el

\textbf{Proof} Let $Y$ be a $\mathbb{G}$-optional process. Suppose that $Y\ind_A$ satisfies the optional splitting formula on the whole time space $\mathbb{R}_+\times\Omega$. Let $Y'=\mathfrak{p}^{[0,\tau)}(Y\ind_A)$ and $Y''=\mathfrak{p}^{[\tau,\infty)}(Y\ind_A)$, respectively. We have$$
Y\ind_A=Y'\ind_{[0,\tau)}+Y''(\tau)\ind_{[\tau,\infty)}.
$$
Since $\ind_A=\ind_A^2$, we also have$$
Y\ind_A=Y\ind_A^2=(Y'\ind_{[0,\tau)}+Y''(\tau)\ind_{[\tau,\infty)})\ind_A,
$$
i.e. the optional splitting formula for $Y$ on $A$.

Conversely, suppose that the optional splitting formula on $A$ holds. Let $C'=\mathfrak{p}^{[0,\tau)}_AY$ and $C''=\mathfrak{p}^{[\tau,\infty)}_AY$, respectively. We then write $$
Y\ind_A=(C'\ind_{[0,\tau)}+C''(\tau)\ind_{[\tau,\infty)})\ind_A.
$$
Note that $A$ is $\mathbb{G}$-predictable. According to Lemma \ref{predictableSPLT}, $\ind_A$ satisfies the global optional splitting formula. Let $B'=\mathfrak{p}^{[0,\tau)}\ind_A$ and $B''=\mathfrak{p}^{[\tau,\infty)}\ind_A$. The above identity becomes$$
\dcb
Y\ind_A
&=&C'B'\ind_{[0,\tau)}+C''(\tau)B''(\tau)\ind_{[\tau,\infty)}.
\dce
$$
This is a global optional splitting formula for $Y\ind_A$. \ok

\bl\label{predictunion}
Let $(A_i)_{i=1}^\infty$ be a sequence of $\mathbb{G}$-predictable sets. Suppose that $(A_i)_{i=1}^\infty\subset\mathcal{L}^o$. Then, $\cup_{i=1}^\infty A_i\in \mathcal{L}^o$.
\el

\textbf{Proof} Let $Y$ be any $\mathbb{G}$-optional process. We apply Lemma \ref{predictLSST} in this proof. It is then enough to prove that $Y\ind_{\cup_{i=1}^{\infty}A_i}$ satisfies the global optional splitting formula. 

By induction, we see that $Y\ind_{\cup_{i=1}^{k}A_i}$ satisfies the global optional splitting formula for any integer $k$. Actually, for $k=1$, this is the case. Suppose that for an integer $k=n$, $Y\ind_{\cup_{i=1}^n A_i}$ satisfies the global optional splitting formula. Let us prove the same for $k=n+1$.
 
We write the identity\!\!:$$
Y\ind_{\cup_{i=1}^{n+1}A_i}=Y\ind_{\cup_{i=1}^{n}A_i}+Y\ind_{A_{n+1}} - Y\ind_{A_{n+1}\cap(\cup_{i=1}^{n}A_i)}.
$$
By assumption, $Y\ind_{\cup_{i=1}^{n}A_i}$ and $Y\ind_{A_{n+1}}$ satisfy the global optional splitting formula. For the term $Y\ind_{A_{n+1}\cap(\cup_{i=1}^{n}A_i)}$, we write it in the form$$
\dcb
Y\ind_{A_{n+1}\cap(\cup_{i=1}^{n}A_i)}
&=&(Y\ind_{\cup_{i=1}^{n}A_i})(\ind_{A_{n+1}}).\\
\dce
$$
$\ind_{A_{n+1}}$ satisfies the global optional splitting formula, because $A_{n+1}$ is $\mathbb{G}$-predictable (cf. Lemma \ref{predictableSPLT}). $Y\ind_{\cup_{i=1}^{n}A_i}$ satisfies the global optional splitting formula by assumption. Applying Lemma \ref{stable}, we conclude that $Y\ind_{\cup_{i=1}^{n+1}A_i}$ also satisfies the global optional splitting formula. 

Now, taking the limit on $Y\ind_{\cup_{i=1}^{k}A_i}$ when $k\rightarrow\infty$, we conclude that $Y\ind_{\cup_{i=1}^{\infty}A_i}$ satisfies the global optional splitting formula (cf. Lemma \ref{stable}). \ok

\subsection{Optional splitting formula on a left-closed right-open interval $[S,T)$}\label{left-closed}

\bl\label{LSST}
Let $S,T$ be two $\mathbb{G}$-stopping times. To have the local optional splitting formula on $[S,T)$, it is necessary and sufficient that, for any bounded $(\mathbb{Q},\mathbb{G})$-martingale $X$ such that $X_T\in\mathcal{G}_{T-}$, $X$ satisfies the optional splitting formula on $[S,T)$.
\el

\textbf{Proof} The condition is necessary by definition. Let us consider the sufficiency. We follow the argument in \cite[Chapitre XX, n$^\circ$22]{DM2}. Let $\xi_\infty\in\mathcal{G}_\infty$ be a bounded random variable and $\xi_t=\mathbb{E}^\mathbb{Q}[\xi|\mathcal{G}_t], t\in\mathbb{R}_+$. Let $$
Y=\Delta_T\xi\ind_{[T,\infty)}-(\Delta_T\xi\ind_{[T,\infty)}){^{\mathbb{G}-(p)}}
$$
and $X=\xi^T-Y$, where $\Delta_T\xi$ denotes the jump of the process $\xi$ at $T$, $\xi^T$ denotes the process $\xi$ stopped at $T$, and $\centerdot{^{\mathbb{G}-(p)}}$ denotes the $(\mathbb{Q},\mathbb{G})$-predictable dual projection. Note that $X,Y$ are $(\mathbb{Q},\mathbb{G})$-martingales. Because $(\Delta_T\xi\ind_{[T,\infty)}){^{\mathbb{G}-(p)}}$ is a $\mathbb{G}$-predictable process, we have$$
\Delta_TX=\Delta_T(\Delta_T\xi\ind_{[T,\infty)}){^{\mathbb{G}-(p)}}\in\mathcal{G}_{T-}
$$
(cf. \cite{Yan}) so that $X_T\in\mathcal{G}_{T-}$. We now write $$
\xi = (\xi-\xi^T)+X+\Delta_T\xi\ind_{[T,\infty)}-(\Delta_T\xi\ind_{[T,\infty)}){^{\mathbb{G}-(p)}}.
$$
For the terms on the right hand side of the above identity, $(\xi-\xi^T)+\Delta_T\xi\ind_{[T,\infty)}$ is null on $[S,T)$ and therefore satisfies the optional splitting formula on $[S,T)$; by assumption $X$ satisfies the optional splitting formula on $[S,T)$; $(\Delta_T\xi\ind_{[T,\infty)}){^{\mathbb{G}-(p)}}$ being $\mathbb{G}$-predictable, satisfies the optional splitting formula thanks to Lemma \ref{predictableSPLT}. Consequently, $\xi$ satisfies the optional splitting formula on $[S,T)$ (cf. Lemma \ref{stable}).

Introduce $\mathfrak{A}$ the family of all bounded $Y\in\mathcal{B}[0,\infty)\otimes\mathcal{G}_\infty$ such that ${^{\mathbb{G}-(o)}}Y$ satisfies the optional splitting formula on $[S,T)$, where ${^{\mathbb{G}-(o)}}\centerdot$ denotes the $(\mathbb{Q},\mathbb{G})$-optional projection.
We can verify that $\mathfrak{A}$ is a functional monotone class (cf. \cite{RW, Yan}) and, according to the above result, $\mathfrak{A}$ contains all the random variables $\ind_{(a,b]}\xi$, where $a,b\in\mathbb{R}_+$ and $\xi$ is a bounded random variable in $\mathcal{G}_\infty$. By the monotone class theorem, $\mathfrak{A}$ contains all bounded $\mathcal{B}[0,\infty)\otimes\mathcal{G}_\infty$-measurable random variables. This implies that all bounded $\mathbb{G}$-optional processes satisfy the local optional splitting formula on $[S,T)$. \ok

\subsection{Local optional splitting formula on the graph of a stopping time}

For a random time $R$, the graph $[R]$ is defined as $[R]=\{(t,\omega)\!\!:t\in\mathbb{R}_+, t=R(\omega)\}$. Note that, if $R=\infty$, $[R]=\emptyset$. By the monotone class theorem we obtain the following lemma.

\bl\label{FR}
Let $R$ be a $\mathbb{G}$-stopping time. For any random variable $\xi\in\mathcal{F}_R$, there exists a $\mathbb{F}$-optional process $Y$ such that $\ind_{\{R<\infty\}}\xi=\ind_{\{R<\infty\}}Y_R$.
\el

In the same vein we have\!\!:

\bl\label{Ytau}
Let $R$ be a $\mathbb{G}$-stopping time. Let $\zeta\in\mathcal{N}\vee\sigma(\tau)\vee\mathcal{F}_\infty$ be a random variable. Then, $\zeta\in\mathcal{N}\vee\sigma(\tau)\vee\mathcal{F}_R$, if and only if there exists a $Y\in\mathcal{B}[0,\infty]\otimes\mathcal{O}(\mathbb{F})$ such that $$
\ind_{\{R<\infty\}}Y_R(\tau)=\ind_{\{R<\infty\}}\zeta.
$$
\el

\textbf{Proof} Let $\hat{\mathcal{C}}$ be the family of all function $Y$ defined on $[0,\infty]\times(\mathbb{R}_+\times\Omega)$ such that $Y_R(\tau)\in\sigma(\tau)\vee\mathcal{F}_R$. $\hat{\mathcal{C}}$ is a functional monotone class, containing the functions $g(t)Z_s(\omega)$, where $g$ is a Borel function on $[0,\infty]$ and $Z\in\mathcal{O}(\mathbb{F})$. By the monotone class theorem, $\hat{\mathcal{C}}$ contains any function $Y$ in $\mathcal{B}[0,\infty]\otimes\mathcal{O}(\mathbb{F})$. 

Suppose the second condition with a $Y\in\mathcal{B}[0,\infty]\otimes\mathcal{O}(\mathbb{F})$. Then, $$
\dcb
\zeta
&=&\ind_{\{R<\infty\}}Y_R(\tau)+\ind_{\{R=\infty\}}\zeta\\
&\in&\{R<\infty\}\cap(\sigma(\tau)\vee\mathcal{F}_R)+\{R=\infty\}\cap(\mathcal{N}\vee\sigma(\tau)\vee\mathcal{F}_R)\\
&\subset&\mathcal{N}\vee\sigma(\tau)\vee\mathcal{F}_R.
\dce
$$

Conversely suppose the first condition. Let $\mathcal{C}$ be the family of all functions on $\Omega$ which satisfy the second condition. Then, $\mathcal{C}$ is a functional monotone class and contains random variables of the form $\ind_Bg(\tau)\xi$, where $B\in\mathcal{N}$, $g$ is a bounded Borel function and $\xi\in\mathcal{F}_R$ (see Lemma \ref{FR}). Applying the monotone class theorem, we conclude that $\mathcal{C}$ contains all $\mathcal{N}\vee\sigma(\tau)\vee\mathcal{F}_R$-measurable random variables. \ok

\bethe\label{graph} 
Let $R$ be a $\mathbb{G}$-stopping time. Then, $[R]\in\mathcal{L}^o$, if and only if $$
\{R<\infty\}\cap\mathcal{G}_R=\{R<\infty\}\cap(\mathcal{N}\vee\sigma(\tau\nmid R)\vee\mathcal{F}_R).
$$
\ethe

\textbf{Proof} Suppose that the local optional splitting formula holds on the graph $[R]$. Let $Y$ be any $\mathbb{G}$-optional process. Let $Y'=\mathfrak{p}^{[0,\tau)}_{[R]}Y$ and $Y''=\mathfrak{p}^{[\tau,\infty)}_{[R]}Y$. Let $\mathcal{T}^0$ be the trivial $\sigma$-algebra\!\!: $\mathcal{T}^0=\{\emptyset,\Omega\}$. Note that $R\in\mathcal{F}_R$ and $\{\tau\leq R<\infty\} \in\{R<\infty\}\cap(\sigma(\tau\nmid R)\vee\mathcal{F}_R)$. We have$$
\dcb
&&Y_R\ind_{\{R<\infty\}}\\
&=&Y'_R\ind_{\{0\leq R<\tau\}}+Y''_R(\tau)\ind_{\{\tau\leq R<\infty\}}\\
&\in&\{0\leq R<\tau\}\cap \mathcal{F}_R+\{\tau\leq R<\infty\}\cap (\sigma(\tau)\vee\mathcal{F}_R)+\{R=\infty\}\cap\mathcal{T}^0\\

&=&\{0\leq R<\tau\}\cap (\sigma(\tau\nmid R)\vee\mathcal{F}_R)+\{\tau\leq R<\infty\}\cap (\sigma(\tau\nmid R)\vee\mathcal{F}_R)+\{R=\infty\}\cap\mathcal{T}^0\\

&=&\{R<\infty\}\cap (\sigma(\tau\nmid R)\vee\mathcal{F}_R)+\{R=\infty\}\cap\mathcal{T}^0.\\
\dce
$$
This measurability relation yields $$
\{R<\infty\}\cap\mathcal{G}_R\subset\{R<\infty\}\cap(\mathcal{N}\vee\sigma(\tau\nmid R)\vee\mathcal{F}_R).
$$
This is actually an equality, because the inverse inclusion is always true. 

Now suppose $\{R<\infty\}\cap\mathcal{G}_R=\{R<\infty\}\cap(\mathcal{N}\vee\sigma(\tau\nmid R)\vee\mathcal{F}_R)$. Let $Y$ be any $\mathbb{G}$-optional process. Since $Y_R\in\mathcal{G}_R$, we have $$
\dcb
\ind_{\{R<\infty\}}Y_R
&\in&\{R<\infty\}\cap(\mathcal{N}\vee\sigma(\tau\nmid R)\vee\mathcal{F}_R)+\{R=\infty\}\cap\mathcal{T}^0\\

&\subset&\{R<\tau\}\cap(\mathcal{N}\vee\mathcal{F}_R)+
\{\tau\leq  R<\infty\}\cap(\mathcal{N}\vee\sigma(\tau)\vee\mathcal{F}_R)+\{R=\infty\}\cap\mathcal{T}^0.
\dce
$$ 
Therefore, there exist $\zeta'\in\mathcal{F}_R$ and $\zeta''\in\sigma(\tau)\vee\mathcal{F}_R$ such that $$
\ind_{\{R<\infty\}}Y_R = \zeta'\ind_{\{R<\tau\}}+\zeta''\ind_{\{\tau\leq R<\infty\}}.
$$
Let $Y'\in\mathcal{O}(\mathbb{F})$ and $Y''\in\mathcal{B}[0,\infty]\otimes\mathcal{O}(\mathbb{F})$ such that $\ind_{\{R<\infty\}}Y'_R=\ind_{\{R<\infty\}}\zeta'$ and $\ind_{\{R<\infty\}}Y''_R(\tau)=\ind_{\{R<\infty\}}\zeta''$ (see Lemma \ref{FR} and the proof of Lemma \ref{Ytau} for the existences of $Y',Y''$). We deduce from the above identity that
$$
\dcb
Y\ind_{[R]}
&=&Y'\ind_{[0,\tau)}\ind_{[R]}+Y''(\tau)\ind_{[\tau,\infty)}\ind_{[R]}.
\dce
$$
$Y$ satisfies the optional splitting formula on $[R]$. \ok

\subsection{Optional splitting formula on intervals such as $[S,T]$ and $(S,T]$}\label{right-closed}

Recall the following notation. Let $T$ be a $\mathbb{G}$ stopping time. Let $A\in\mathcal{G}_T$. We denote $T_{A}=T\ind_A +\infty\ind_A$ (called the restriction of $T$ on $A$). $T_A$ is again a $\mathbb{G}$ stopping time (cf. \cite{Yan}).

\bl\label{interval+pointLSST++}
Let $S,T$ be $\mathbb{G}$-stopping times. Suppose that $[S,T)\in\mathcal{L}^o$ and $[T_{\{S\leq T<\infty\}}]\in\mathcal{L}^o$. Suppose that $\ind_{[T_{\{S\leq T<\infty\}}]}$ satisfies the optional splitting formula on $[S,T]$. Then, $[S,T]\in\mathcal{L}^o$. 
\el

\textbf{Proof} Let $Y$ be a $\mathbb{G}$-optional process. Let $A'=\mathfrak{p}^{[0,\tau)}_{[S,T)}Y$ and $A''=\mathfrak{p}^{[\tau,\infty)}_{[S,T)}Y$. Let $B'=\mathfrak{p}^{[0,\tau)}_{[T_{\{S\leq T<\infty\}}]}Y$ and $B''=\mathfrak{p}^{[\tau,\infty)}_{[T_{\{S\leq T<\infty\}}]}Y$. Let $C'=\mathfrak{p}^{[0,\tau)}_{[S,T]}\ind_{[T_{\{S\leq T<\infty\}}]}$ and $C''=\mathfrak{p}^{[\tau,\infty)}_{[S,T]}\ind_{[T_{\{S\leq T<\infty\}}]}$. Note that $\ind_{[T_{\{S\leq T<\infty\}}]}=\ind_{[T]}\ind_{\{S\leq T<\infty\}}$. We can write $$
\dcb
Y\ind_{[S,T]}
&=&Y\ind_{[S,T)}+Y\ind_{[T]}\ind_{\{S\leq T<\infty\}}\\

&=&(A'\ind_{[0,\tau)}+A''\ind_{[\tau,\infty)})\ind_{[S,T)}+(B'\ind_{[0,\tau)}+B''\ind_{[\tau,\infty)})\ind_{[T_{\{S\leq T<\infty\}}]}\\

&=&(A'\ind_{[0,\tau)}+A''\ind_{[\tau,\infty)})\ind_{[S,T]}(1-\ind_{[T_{\{S\leq T<\infty\}}]})+(B'\ind_{[0,\tau)}+B''\ind_{[\tau,\infty)})\ind_{[T_{\{S\leq T<\infty\}}]}\ind_{[S,T]}\\

&=&((A'+(B'-A')C')\ind_{[0,\tau)}+(A''+(B''-A'')C'')\ind_{[\tau,\infty)})\ind_{[S,T]}. \ \ \ok
\dce
$$

In the same way we can prove

\bl\label{interval+pointLSST}
Let $S,T$ be $\mathbb{G}$-stopping times. Suppose that $(S,T)\in\mathcal{L}^o$ and $[T_{\{S< T<\infty\}}]\in\mathcal{L}^o$.  Suppose that $\ind_{[T_{\{S<T<\infty\}}]}$ satisfies the optional splitting formula on $(S,T]$. Then, $(S,T]\in\mathcal{L}^o$. 
\el

\

\section{Sufficient conditions to have optional splitting formulas at $\tau$ with respect to $\mathbb{F}$}\label{suff-cond}

\subsection{Optional splitting formula on $[0,\tau)$}

\bethe\label{beforedefault}
$[0,\tau)\in\mathcal{L}^o$.
\ethe

\textbf{Proof} Let $\xi_\infty\in\mathcal{G}_\infty$ be a bounded random variable and $\xi_t=\mathbb{E}^\mathbb{Q}[\xi|\mathcal{G}_t], t\in\mathbb{R}_+$. We write the identity (cf. \cite{DM2,JR})\!\!:
$$
\xi_t\ind_{\{t<\tau\}}
=\ind_{\{t<\tau\}}\frac{\mathbb{Q}[\xi\ind_{\{t<\tau\}}|\mathcal{F}_t]}{\mathbb{Q}[t<\tau|\mathcal{F}_t]}\ind_{\{\mathbb{Q}[t<\tau|\mathcal{F}_t]>0\}}, t\geq 0.
$$
This is an optional splitting formula for $\xi$ on $[0,\tau)$. Now, applying Lemma \ref{LSST}, we conclude the theorem. \ok

From this theorem we deduce the result.

\bcor\label{R<tau}
Let $R$ be a $\mathbb{G}$-stopping time. We have
$$
\{R<\tau\}\cap\mathcal{G}_R = \{R< \tau\}\cap(\mathcal{N}\vee\mathcal{F}_R)
$$
\ecor

\textbf{Proof.}
We have the identity$$
\ind_{[R]}\ind_{[0,\tau)}
=\ind_{[R]}\ind_{\{R<\tau\}}
=\ind_{[R_{\{R<\tau\}}]}.
$$
Since $[0,\tau)\in\mathcal{L}^o$ by Theorem \ref{beforedefault}, $[R_{\{R<\tau\}}]\in\mathcal{L}^o$ (Lemma \ref{on-subset}). According to Theorem \ref{graph}, $$
\{R_{\{R<\tau\}}<\infty\}\cap\mathcal{G}_{R_{\{R<\tau\}}}=\{R_{\{R<\tau\}}<\infty\}\cap(\mathcal{N}\vee\sigma(\tau\nmid R_{\{R<\tau\}})\vee\mathcal{F}_{R_{\{R<\tau\}}}),
$$
which is equivalent to
$
\{R<\tau\}\cap\mathcal{G}_{R}=\{R<\tau\}\cap(\mathcal{N}\vee\mathcal{F}_{R})
$. \ok

\subsection{$s\!\mathcal{H}$-measure}

\begin{definition}\label{Hmeasure}
Let $S,T$ be $\mathbb{G}$-stopping times. A probability measure $\mathbb{Q}'$ defined on $\mathcal{G}_\infty$ is called an $s\!\mathcal{H}$-measure over the random time interval $(S,T]$ (with respect to $(\mathbb{Q},\mathbb{F},\mathbb{G})$), if $\mathbb{Q}'$ is equivalent to $\mathbb{Q}$ on $\mathcal{G}_{\infty}$, and if, for any $(\mathbb{Q},\mathbb{F})$ local martingale $X$, $X^{(S,T]}$ is a $(\mathbb{Q}',\mathbb{G})$ local martingale, where $X^{(S,T]}_t=X^{S\vee T}_t-X^S_t, t\geq 0$. 
\end{definition}

\brem\label{rmkdef}
The notion of $s\!\mathcal{H}$-measure is derived from the general study of the enlargement of filtration in \cite{SongThesis, song-local-solution}. It is employed in \cite{JS3} to study the martingale representation property in $\mathbb{G}$. The above Definition \ref{Hmeasure} is a different but equivalent version of that used in \cite{JS3} (see Lemma A.5 and Lemma A.6 in \cite{JS3}).

Note also that, if $\mathbb{Q}'$ is an $s\!\mathcal{H}$-measure on $(S,T]$ and if $(S',T']\subset(S,T]$, then $\mathbb{Q}'$ is an $s\!\mathcal{H}$-measure on $(S',T']$.
\erem

\brem\label{invariant}
Note that the property of the optional splitting formula on a $\mathbb{G}$-optional set is invariant by the equivalent changes of probability measures on $\mathcal{G}_{\infty}$.
\erem

\bethe\label{gSTMRT}
For any $\mathbb{F}$-stopping time  $T$, for any $\mathbb{G}$-stopping time $S$ such that $S\geq \tau$ (an almost sure relation), if an $s\!\mathcal{H}$-measure $\mathbb{Q}'$ over $(S,T]$ exists, then $[S,T)\in\mathcal{L}^o$.
\ethe
\noindent\textbf{Proof.} Let $\mathbb{Q}'$ be an $s\!\mathcal{H}$-measure over $(S,T]$. Let $\zeta$ be a $\mathcal{F}_{T}$-measurable bounded random variable. We introduce the martingale $X_t=\mathbb{E}^\mathbb{Q}[\zeta|\mathcal{F}_t], t\geq 0$. We note that $\zeta=X_t$ for all $t\geq T$.

Since $X$ is bounded, $X^{(S,T]}$ is a bounded $(\mathbb{Q}',\mathbb{G})$ martingale. Hence,$$
\mathbb{Q}'[X^{(S,T]}_\infty|\mathcal{G}_t]=X^{(S,T]}_t=X^{S\vee T}_{t}-X^S_t, t\geq 0.
$$
The relation $\tau\leq S$ (which implies $\sigma(\tau)\vee\mathcal{F}_S\subset\mathcal{G}_{S}$), which holds under $\mathbb{Q}$, remains valid under $\mathbb{Q}'$. Let $g$ be a bounded Borel function. Noting that $g(\tau)\ind_{\{S<t\}}\in\mathcal{G}_t$, we have$$
\dcb
\mathbb{Q}'[g(\tau)X_{S\vee T}-g(\tau)X_S|\mathcal{G}_t]
&=&\mathbb{Q}'[g(\tau)X^{(S,T]}_\infty|\mathcal{G}_{S\vee t}|\mathcal{G}_t]\\

&=&\mathbb{Q}'[g(\tau)X^{(S,T]}_{t}\ind_{\{S<t\}}|\mathcal{G}_t]\\
&=&g(\tau)X^{S\vee T}_{t}-g(\tau)X^S_t, t\geq 0.
\dce
$$
Consider this identity on the set $\{S\leq t<T\}$. We obtain$$
\ind_{\{S\leq t<T\}}g(\tau)X_t=\ind_{\{S\leq t<T\}}\mathbb{Q}'[g(\tau)\zeta|\mathcal{G}_t].
$$
This identity means that the $(\mathbb{Q}',\mathbb{G})$-martingale $\mathbb{Q}'[g(\tau)\zeta|\mathcal{G}_t], t\geq 0,$ satisfies the optional splitting formula on $[S,T)$. 

Let $\mathcal{C}$ be the class of all bounded random variables $\xi\in\mathcal{G}_\infty$ such that the martingale $\mathbb{Q}'[\xi|\mathcal{G}_t], t\geq 0,$ satisfies the optional splitting formula on $[S,T)$. The preceding result, together with the monotone class theorem, implies that $\mathcal{C}$ contains all bounded $\sigma(\tau)\vee \mathcal{F}_{T}$ measurable random variables. By \cite[Lemme(4.4)]{J}, $
\mathcal{G}_{T -}\subset (\mathcal{N}\vee\sigma(\tau)\vee \mathcal{F}_{T})
$
(noting that we have the same family of negligible sets under $\mathbb{Q}$ and under $\mathbb{Q}'$).
Lemma \ref{LSST} is applicable to conclude the optional splitting formula on $[S,T)$ under the probability measure $\mathbb{Q}'$. Finally, Remark \ref{invariant} completes the proof. \ok

\subsection{Structure of $\mathcal{G}_R$ under an $s\!\mathcal{H}$-measure}

\bl\label{pointLSST}
Let $R$ be a $\mathbb{G}$-stopping time. For any $\mathbb{F}$-stopping time  $T$ and any $\mathbb{G}$-stopping time $S$, if an $s\!\mathcal{H}$-measure $\mathbb{Q}'$ over $(S,T]$ exists, we have$$
\{\tau\leq R\}\cap\{S\leq R<T\}\cap\mathcal{G}_R = \{\tau \leq R\}\cap\{S\leq R<T\}\cap(\mathcal{N}\vee\sigma(\tau)\vee\mathcal{F}_R).
$$
 
\el

\textbf{Proof} Let $\mathbb{Q}'$ be an $s\!\mathcal{H}$-measure over $(S,T]$. Let $\zeta$ be a $\mathcal{F}_{T}$-measurable bounded random variable. We introduce the martingale $X_t=\mathbb{E}^\mathbb{Q}[\xi|\mathcal{F}_t], t\geq 0$. Let $g$ be a bounded Borel function. As in the previous lemma, we prove
$$
\dcb
\ind_{\{\tau\leq R\}}\ind_{\{S\leq R<T\}}g(\tau)X_{R}

&=&\ind_{\{\tau\leq R\}}\ind_{\{S\leq R<T\}}\mathbb{Q}'[g(\tau)\zeta|\mathcal{G}_{R}].
\dce
$$
Then we can write $$
\dcb
&&\ind_{\{\tau\leq R\}}\ind_{\{S\leq R<T\}}\mathbb{Q}'[g(\tau)\zeta|\mathcal{G}_{R}]\\
&\in&\{\tau\leq R\}\cap\{S\leq R<T\}\cap(\sigma(\tau)\vee\mathcal{F}_{R})
+(\{\tau\leq R\}\cap\{S\leq R<T\})^c\cap\mathcal{T}^0.
\dce
$$
($\mathcal{T}^0$ denotes the trivial $\sigma$-algebra.) By the monotone class theorem, this relation is extended to any bounded $\xi\in\sigma(\tau)\vee\mathcal{F}_{T}$\!\!:
\begin{equation}\label{incl1}
\dcb
&&\ind_{\{\tau\leq R\}}\ind_{\{S\leq R<T\}}\mathbb{Q}'[\xi|\mathcal{G}_{R}]\\
&\in&\{\tau\leq R\}\cap\{S\leq R<T\}\cap(\sigma(\tau)\vee\mathcal{F}_{R})
+(\{\tau\leq R\}\cap\{S\leq R<T\})^c\cap\mathcal{T}^0.
\dce
\end{equation}
We note that $\{\tau\leq R\}\cap\{S\leq R<T\}\in\mathcal{G}_{T-}$ and $$
\dcb
\{\tau\leq R\}\cap\{S\leq R<T\}\cap\mathcal{G}_R&\subset&
\{\tau\leq R\}\cap\{S\leq R<T\}\cap\mathcal{G}_{T-}\\
&\subset&\mathcal{N}\vee\sigma(\tau)\vee\mathcal{F}_{T-},
\dce
$$
where the last inclusion is a consequence of \cite[Lemme(4.4)]{J}. Applying the equation (\ref{incl1}) to all the $\xi\in\mathcal{G}_R$, we conclude that $$
\{\tau\leq R\}\cap\{S\leq R<T\}\cap\mathcal{G}_R\subset
\{\tau\leq R\}\cap\{S\leq R<T\}\cap(\mathcal{N}\vee\sigma(\tau)\vee\mathcal{F}_{R}).
$$
The inverse inclusion relation being obvious, we actually have an equality. \ok

\subsection{$s\!\mathcal{H}$-measure with covering condition}

\bethe\label{mrt_after_default}
Suppose that there exists a countable family of pairs of $\mathbb{G}$ stopping times $\{S_j,T_j\}, j\in \mathbb{N},$ such that
\ebe
\item[(1)]
$T_j$ are $\mathbb{F}$-stopping times;
 
\item[(2)]
$(\tau,\infty)\subset\cup_{i\in\mathbb{N}}(S_j,T_j)$ (covering condition on $(\tau,\infty)$).

\dbe
Suppose that, for any $j\in\mathbb{N}$, there exists an $s\!\mathcal{H}$-measure $\mathbb{Q}_j$ over the time interval $(S_j,T_j]$. Then $(\tau,\infty)\in\mathcal{L}^o$.

If we replace the condition $(2)$ with the condition\!\!:
\ebe
\item[(2)']
$[\tau,\infty)\cap(0,\infty)\subset\cup_{i\in\mathbb{N}}(S_j,T_j)$ (the covering condition on $[\tau,\infty)$),
\dbe
then the global optional splitting formula holds. 
\ethe

\textbf{Proof} Let us suppose the covering condition $(2)$. For any $k\geq 0$, consider the pair $\{S_k,T_k\}$. Let us verify the three conditions in Lemma \ref{interval+pointLSST} with respect to the random interval $(\tau\vee S_k,T_k]$. First, according to Theorem \ref{gSTMRT}, $[\tau\vee S_k,T_k)\in\mathcal{L}^o$, and \textit{à fortiori} $(\tau\vee S_k,T_k)\in\mathcal{L}^o$. Second, we apply Lemma \ref{pointLSST} to write the equality\!\!: $j\geq 1$,$$
\{\tau\leq T_k\}\cap\{S_j\leq T_k<T_j\}\cap\mathcal{G}_{T_k} = \{\tau \leq T_k\}\cap\{S_j\leq T_k<T_j\}\cap(\mathcal{N}\vee\sigma(\tau)\vee\mathcal{F}_{T_k}).
$$
From this, we get
$$
\{\tau< T_k\}\cap\{S_j< T_k<T_j\}\cap\mathcal{G}_{T_k} = \{\tau < T_k\}\cap\{S_j< T_k<T_j\}\cap(\mathcal{N}\vee\sigma(\tau)\vee\mathcal{F}_{T_k}).
$$
Since we have $$
\dcb
\{S_j<T_k<T_j\}\in\mathcal{G}_{T_k},\\
\{T_k<T_j\}\in\mathcal{F}_{T_k}\subset \sigma(\tau)\vee\mathcal{F}_{T_k},\\
\{S_j<T_k\}\in\mathcal{G}_{T_k-}\subset \mathcal{N}\vee\sigma(\tau)\vee\mathcal{F}_{T_k},
\dce
$$
Lemma \ref{union-AT} is applicable. By then taking the union on $j\geq 0$, using the covering condition $(2)$, we obtain $$
\{\tau< T_k<\infty\}\cap\mathcal{G}_{T_k} = \{\tau < T_k<\infty\}\cap(\mathcal{N}\vee\sigma(\tau)\vee\mathcal{F}_{T_k}).
$$
By Theorem \ref{graph}, $[(T_k)_{\{\tau< T_k\}}]\in\mathcal{L}^o$, and \textit{à fortiori} $[(T_k)_{\{\tau\vee S_k< T_k<\infty\}}]\in\mathcal{L}^o$. Finally, we write $$
\ind_{[(T_k)_{\{\tau\vee S_k< T_k<\infty\}}]}
=\ind_{[T_k]}\ind_{\{\tau\vee S_k< T_k<\infty\}}
=\ind_{[T_k]}\ind_{(\tau\vee S_k,T_k]}.
$$
Since $\ind_{[T_k]}$ is a $\mathbb{F}$ optional process, since $(\tau\vee S_k,T_k]\subset[\tau,\infty)$, the above formula is an optional splitting formula for $\ind_{[(T_k)_{\{\tau\vee S_k< T_k<\infty\}}]}$ on $(\tau\vee S_k,T_k]$. 

According to Lemma \ref{interval+pointLSST}, $(\tau\vee S_k,T_k]\in\mathcal{L}^o$ for $k\geq 0$. Also they are $\mathbb{G}$-predictable sets. By Lemma \ref{predictunion} and the covering condition (2), $(\tau,\infty)=\cup_{k\in\mathbb{N}}(\tau\vee S_k,T_k]\in\mathcal{L}^o$. The first part of the theorem is proved.

Now suppose the covering condition $(2)'$. Applying Lemma \ref{pointLSST} to the random time $\tau$, we write$$
\{S_j< \tau<T_j\}\cap\mathcal{G}_{\tau} = \{S_j< \tau<T_j\}\cap(\mathcal{N}\vee\sigma(\tau)\vee\mathcal{F}_{\tau}).
$$
Since $\{S_j< \tau<T_j\}\in\mathcal{G}_{\tau}$, $\{S_j< \tau\}\in \mathcal{G}_{\tau-}\subset \sigma(\tau)\vee\mathcal{F}_{\tau}$, and $\{\tau<T_j\}\in\mathcal{F}_{\tau}$, Lemma \ref{union-AT} is applicable. Taking the union on $j\geq 0$, using the covering condition $(2)'$, we obtain $$
\{0< \tau<\infty\}\cap\mathcal{G}_{\tau} 
= \{0< \tau<\infty\}\cap(\mathcal{N}\vee\sigma(\tau)\vee\mathcal{F}_{\tau})
=\{0< \tau<\infty\}\cap(\mathcal{N}\vee\mathcal{F}_{\tau}).
$$
On the other side, according to \cite[Lemme(4.4)]{J},$$
\{\tau=0\}\cap\mathcal{G}_{\tau}
=\{\tau=0\}\cap\mathcal{G}_{0}
=\{\tau=0\}\cap(\mathcal{N}\vee\mathcal{F}_{0}) 
=\{\tau=0\}\cap(\mathcal{N}\vee\mathcal{F}_{\tau}).
$$
Taking the union of these two identities, we obtain
$$
\{\tau<\infty\}\cap\mathcal{G}_{\tau} 
=\{\tau<\infty\}\cap(\mathcal{N}\vee\mathcal{F}_{\tau}).
$$
Theorem \ref{graph} is now applicable to conclude $[\tau]\in\mathcal{L}^o$. 

Let $Y$ be a $\mathbb{G}$ optional process. Let $A''=\mathfrak{p}^{[\tau,\infty)}_{(\tau,\infty)}Y$ and $B''=\mathfrak{p}^{[\tau,\infty)}_{[\tau]}Y$. We check right away that$$
\dcb
Y\ind_{[\tau,\infty)}
&=&(A''\ind_{(\tau,\infty)}+B''\ind_{[\tau]})\ind_{[\tau,\infty)}.
\dce
$$
This proves the local splitting formula on $[\tau,\infty)$. Since the local splitting formula always holds on $[0,\tau)$ (Theorem \ref{beforedefault}), the second part of the theorem is proved. \ok

\

\section{Examples}
\label{examples}

In this section we expose the connection between the results in the previous sections and different works in the literature of credit risk modeling.

\subsection{Hypothesis$(H)$}\label{dissH}

We say that the hypothesis$(H)$ is satisfied between the pair of filtrations $(\mathbb{F},\mathbb{G})$ under $\mathbb{Q}$, if every $(\mathbb{Q},\mathbb{F})$-martingale is a $(\mathbb{Q},\mathbb{G})$-martingale. This hypothesis$(H)$ has been used in numerous papers on credit risk modeling. The hypothesis$(H)$ can be characterized with different equivalent conditions (cf. \cite{BJR,BY}). In particular, the hypothesis$(H)$ is satisfied if $\tau$ is independent of $\mathcal{F}_\infty$ or if $\tau$ is a Cox time (cf. \cite{BJR}).

\bethe\label{HyH}
If there exists a probability measure $\mathbb{Q}'$ equivalent to $\mathbb{Q}$ such that the hypothesis$(H)$ is satisfied under $\mathbb{Q}'$, then the probability measure $\mathbb{Q}'$ is an $s\!\mathcal{H}$-measure over $(0,\infty]$. Consequently, the global optional splitting formula holds. 
\ethe

\textbf{Proof.} The first part of the theorem can be checked by definition. For the second part, we note that the covering condition $(2')$ is satisfied. The global optional splitting formula, therefore, is the consequence of Theorem \ref{mrt_after_default}. \ok

Note that the condition of the above theorem is satisfied in the case of the hypothesis$(H)$ or in the proportionality model \cite{JS}. The theorem also is applicable in a model satisfying the density hypothesis (cf. Sections \ref{multi-time} and \ref{multi-time-marked} below), if its density function is strictly positive. This said, better results can be proved on the density hypothesis. See Lemma \ref{one-tau-density} and Lemma \ref{one-tau-density-marked}. See \cite{CJZ, ElKJJ, JL, jiao, KL, Pham} for applications under the density hypothesis.

\subsection{Review of some results}

We now return to the works of \cite{BSW, biagini, Kusuoka} mentioned in subsection \ref{bgd} and show that these results can be proved with the optional splitting formula. 

The work \cite{BSW} established the right-continuity of the filtration of $\sigma(\tau\nmid t)\vee\mathcal{F}_t$ (completed by the null sets), $t\geq 0$, when $\tau$ is a Cox time. For a new proof, we know from the previous subsection \ref{dissH} that a Cox time satisfies hypothesis$(H)$. According to Theorem \ref{HyH} the global optional splitting formula holds. Applying Theorem \ref{right-continuity} we obtain $\mathcal{G}_t=\mathcal{N}\vee\sigma(\tau\nmid t)\vee\mathcal{F}_t, t\geq 0$. The result of \cite{BSW} is proved, because $\mathbb{G}$ is a right-continuous filtration.

In the proof of \cite[Proposition 4.1]{biagini} it is found that, for any $\mathbb{G}$-martingale $Z$, there exists a $\mathbb{F}$-predictable process $\hat{Z}$ such that $Z_\tau=\hat{Z}_\tau$ if $\tau<\infty$. This is equivalent to saying that $\{\tau<\infty\}\cap\mathcal{G}_{\tau}=\{\tau<\infty\}\cap\mathcal{G}_{\tau-}$. Let us explain this property with the optional splitting formula. Indeed, the assumption of \cite{biagini} implies that the hypothesis($H$) holds as well as the global optional splitting formula. Hence, $[\tau]\in\mathcal{L}^o$, according to Theorem \ref{graph}, $$
\{\tau<\infty\}\cap\mathcal{G}_{\tau}=\{\tau<\infty\}\cap(\mathcal{N}\vee\sigma(\tau)\vee\mathcal{F}_{\tau}).
$$
If $\mathbb{F}$ is moreover a Brownian filtration as assumed in \cite{biagini}, $$
\mathcal{N}\vee\sigma(\tau)\vee\mathcal{F}_{\tau}
=\mathcal{N}\vee\sigma(\tau)\vee\mathcal{F}_{\tau-}
=\mathcal{G}_{\tau-}.
$$
We obtain the desired equality.

In \cite{Kusuoka} a random time $\tau$ is considered, with a continuous probability distribution and satisfying the hypothesis$(H)$. It is then assumed that the filtration $\sigma(\tau\wedge t)\vee\mathcal{F}_t, t\geq 0$, is right-continuous. Let us show that, under the assumptions of \cite{Kusuoka}, there is no need to assume this right-continuity, because it is automatically true. Actually, since $\tau$ has a continuous distribution, $$
\mathcal{N}\vee\sigma(\tau\wedge t)\vee\mathcal{F}_t
=\mathcal{N}\vee\sigma(\tau\nmid t)\vee\mathcal{F}_t.
$$
Now applying Theorem \ref{right-continuity} (passing through Theorem \ref{HyH}), $\mathcal{N}\vee\sigma(\tau\nmid t)\vee\mathcal{F}_t$ coincides with $\mathcal{G}_t$, which is right-continuous.

\subsection{Honest time} 

A random time is called honest if it is equal to the end of an optional set, when it is finite. A large literature exists on the subject of honest time. We mention, among many others, some of the first papers \cite{barlow, J, JY} and some of the applications in financial modeling \cite{FJS,NP}. 

Honest time was the first (counter-)example in which the problem with the global optional splitting formula was revealed. We can mention the example in \cite{barlow} (see subsection \ref{subject} for a description of the example). We also mention \cite[Proposition(5.6)]{J} (see subsection \ref{some-consequences} for a description) which generalizes the example in \cite{barlow}. Notice, however, according to \cite[Proposition(5.3) b)]{J}, $(\tau,\infty)\in\mathcal{L}^o$ for a honest time $\tau$.

\subsection{$\natural$-model}

We present in this subsection a model developed in \cite{JS2}. Through this example we explain how to check the optional splitting formula when the enlargement of filtration formula is known. 

It is a model on a product probability space. We are given a space $\hat{\Omega}$ equipped with a filtration $\hat{\mathbb{F}}=(\hat{\mathcal{F}}_t)_{t\in\mathbb{R}_+}$ and a probability measure $\hat{\mathbb{Q}}$ on $\hat{\mathcal{F}}_\infty$. We consider the product space $\Omega=[0,\infty]\times\hat{\Omega}$ equipped with the product $\sigma$-algebra $\mathcal{B}[0,\infty]\otimes\hat{\mathcal{F}}_\infty$. Let $\pi$ and $\tau$ denote the projection maps \!: $\pi(s,\hat{\omega})=\hat{\omega}$ and $\tau(s,\hat{\omega})=s$ for $(s,\hat{\omega})\in\Omega$. We define $\mathbb{F}$ to be the filtration $\mathcal{F}_t=\pi^{-1}(\hat{\mathcal{F}}_t), t\in\mathbb{R}_+$, and $\mathbb{Q}$ to be the probability measure on $\mathcal{F}_\infty$ such that $\mathbb{Q}(\pi^{-1}(\mathtt{A}))=\hat{\mathbb{Q}}(\mathtt{A})$ for $\mathtt{A}\in\hat{\mathcal{F}}_\infty$. The triplet $(\mathbb{Q},\mathbb{F},\tau)$ represents a financial market with a credit default time $\tau$. 

The problem considered in \cite{JS2} is the following. We are given an $\mathbb{F}$-adapted continuous increasing process $\Lambda$ and a non negative $(\mathbb{Q},\mathbb{F})$ local martingale $N$, such that $\Lambda_0=0,N_0=1$ and $0<N_t e^{-\Lambda_t}< 1$ for all $t\in\mathbb{R}_+$. 

\textbf{Problem $\mathcal{P}^*$.} Construct on $\mathcal{B}[0,\infty]\otimes\hat{\mathcal{F}}_\infty$ a probability measure $\widetilde{\mathbb{Q}}$ such that 
\ebe
\item[$-$]
(restriction condition) $\widetilde{\mathbb{Q}}|_{\mathcal{F}_\infty}=\mathbb{Q}|_{\mathcal{F}_\infty}$ and
\item[$-$]
(projection condition) $\widetilde{\mathbb{Q}}[\tau>t|\mathcal{F}_t]=N_t
e^{-\Lambda_t}$ for all $t\in\mathbb{R}_+$.
\dbe
The problem $\mathcal{P}^*$ is essential for a useful theory of credit default modeling through the progressive enlargement of filtration. The process $N_t
e^{-\Lambda_t}, t\geq 0,$ represents the data calibrated from the market. We need to know whether, for any type of market, there exists a corresponding credit default model. 

For long time, the problem $\mathcal{P}^*$ was only solved in the case where $N\equiv 1$ via Cox process method (cf. \cite{BJR}). In \cite{JS2} the following result is proved. We suppose that all $(\mathbb{Q},\mathbb{F})$ local martingales are continuous. 
Then, for any $(\mathbb{Q},\mathbb{F})$ local martingale $Y$, for any bounded differentiable function $f$ with bounded continuous derivative and $f(0)=0$, there exists $\mathbb{Q}^\natural$ solving the problem $\mathcal{P}^*$ on the product space such that, for any $u\in\mathbb{R}_+^*$, the martingale $M^u_t=\mathbb{Q}^\natural[\tau\leq u|\mathcal{F}_t], t\geq u,$ satisfies the following evolution equation($\natural$)\!:$$
(\natural_u)  \left\{\dcb
dX_t=X_t\left(-\frac{e^{-\Lambda_t}}{1-Z_t}dN_t+f(X_t - (1-Z_t))dY_t\right),\ t\in[u,\infty),\\
X_u=1-Z_u,
\dce
\right.
$$
where $Z_t=N_te^{-\Lambda_t}$. 

As a consequence, there exists an infinity of solutions to the problem $\mathcal{P}^*$. Moreover, if in addition, for $0< t< \infty$, the map $u\rightarrow M^u_t$ is continuous on $(0,t]$, then any $(\mathbb{Q},\mathbb{F})$ local martingale $X$ is a $(\mathbb{Q}^\natural,\mathbb{G})$ semimartingale in the way that $\widetilde{X}=X-\Gamma(X)$ is a $(\mathbb{Q}^\natural,\mathbb{G})$ local martingale, where $\Gamma(X)$ (the drift process) is given by the following formula (called enlargement of filtration formula):
\begin{equation}\label{decomposition-formula}
\dcb
\pv(X)_t=\int_{0}^{t}\ind_{\{s\leq \tau\}}\frac{e^{-\Lambda_s}}{Z_s}d\cro{N,X}_s
-\int_{0}^{t}\ind_{\{\tau< s\}}\frac{e^{-\Lambda_s}}{1-Z_s}d\cro{N,X}_s\\
\\
\hspace{1.5cm}+\int_{0}^{t}\ind_{\{\tau< s\}}( f(M^{\tau}_s-(1-Z_s))+M^\tau_s f'(M^{\tau}_s-(1-Z_s))) d\cro{Y,X}_s,\ t\in\mathbb{R}_+.
\dce
\end{equation}
It it interesting to note that, if $f\equiv 0$, the above formula takes exactly the same form than that in the case of a honest time (cf. \cite{J}). However, unlike the honest time model, the optional splitting formula holds in this $\natural$-model.  

Following Theorem \ref{mrt_after_default} we look for $s\!\mathcal{H}$-measures. Notice that an $s\!\mathcal{H}$-measure is simply a probability change which cancels locally the drift process $\Gamma(X)$ in the filtration $\mathbb{G}$. We introduce $$
\gamma_s=\frac{e^{-\Lambda_s}}{Z_s},\
\alpha_s=-\frac{e^{-\Lambda_s}}{1-Z_s},\
\beta_s=f(M^{\tau}_s-(1-Z_s))+M^\tau_s f'(M^{\tau}_s-(1-Z_s)).
$$
With these notations we can write the drift process $\Gamma(X)$ in the form$$
d\Gamma(X)_t=
(\gamma_t\ind_{\{t\leq \tau\}}+\alpha_t\ind_{\{\tau<t\}})d\cro{N,X}
+\beta_t\ind_{\{\tau<t\}}d\cro{Y,X}.
$$
Recall the notations $\widetilde{X}=X-\Gamma(X), \widetilde{N}=N-\Gamma(N),\widetilde{Y}=Y-\Gamma(Y)$ which are $(\mathbb{Q}^\natural,\mathbb{G})$ local martingales. By the continuity, we have $$
\cro{N,X}=\cro{\widetilde{N},\widetilde{X}},\ \cro{Y,X}=\cro{\widetilde{Y},\widetilde{X}},
$$ 
so that
$$
d\Gamma(X)_t=
(\gamma_t\ind_{\{t\leq \tau\}}+\alpha_t\ind_{\{\tau<t\}})d\cro{\widetilde{N},\widetilde{X}}
+\beta_t\ind_{\{\tau<t\}}d\cro{\widetilde{Y},\widetilde{X}}.
$$
This expression of $\Gamma(X)$ in term of $\mathbb{G}$ local martingales indicates how to use Girsanov's theorem to cancel locally the drift process. Therefore, for $0<a<\infty,n\in\mathbb{N}^\ast$, we introduce $$
\dcb
T_{a,n}&=\inf\{v\geq a\!\!:\hspace{3pt} & \int_a^v\gamma^2_s+\alpha^2_sd\cro{N}_w+\int_a^vd\cro{Y}+(v-a)>n\},
\dce
$$
which is a $\mathbb{F}$-stopping time, and we define the exponential martingale $$
\eta^{a,n}_t=\mathcal{E}\left(\int_0^t(-\gamma_s\ind_{\{s\leq \tau\}}-\alpha_s\ind_{\{\tau<s\}})\ind_{\{a<s\leq T_{a,n}\}}d\widetilde{N}_s
+\int_0^t(-\beta_s)\ind_{\{\tau<s\}}\ind_{\{a<s\leq T_{a,n}\}} d\widetilde{Y}_s\right),
$$
$t\in\mathbb{R}_+$, whose associated probability measure, making $X$ a $(\mathbb{Q}^\natural,\mathbb{G})$ local martingale by Girsanov's formula, is an $s\!\mathcal{H}$-measure on $(a, T_{a,n}]$. 

Notice that usually we would define the stopping times $T_{a,n}$ with the process $(-\gamma\ind_{[0,\tau]}-\alpha\ind_{(\tau,\infty)})^2$. But we can not do so, because Theorem \ref{mrt_after_default} requires $T_{a,n}$ to be $\mathbb{F}$ stopping time. Since $Z,N,Y$ are continuous and $0<Z<1$ on $(0,\infty)$, $\lim_{n\rightarrow\infty}T_{a,n}=\infty$ which implies $(0,\infty)=\cup_{a\in\mathtt{Q},n\in\mathbb{N}^*}(a,T_{a,n})$. The $s\!\mathcal{H}$-measure condition covering $[\tau,\infty)$ in Theorem \ref{mrt_after_default} is satisfied. Consequently, the global optional splitting formula at the random time $\tau$ holds in this $\natural$-model. We emphasize that in this example the hypothesis$(H)$ is not involved.

\

\section{Splitting formula at multiple random times}\label{multi-time}

In this section we extend the preceding results to the case of multiple random times.

\subsection{Ordering the values of a positive function defined on the set $\{1,\ldots,k\}$}\label{ordering}

Let $\mathfrak{a}$ be a function defined on $\{1,\ldots,k\}$  (where $k>0$ is an integer) taking values in $[0,\infty]$. Let $\{a_1,\ldots,a_k\}$ denote the values of $\mathfrak{a}$. Consider the points $(a_1,1),\ldots,(a_k,k)$ in the space $[0,\infty]\times\{1,\ldots,k\}$. These points are two-by-two distinct. We order these points according to alphabetic order in the space $[0,\infty]\times\{1,\ldots,k\}$. Then, for $1\leq i\leq k,$ the rang of $(a_i,i)$ in this ordering is given by $$
R^\mathfrak{a}(i)=R^{\{a_1,\ldots,a_k\}}(i)=\sum_{j=1}^k \ind_{\{a_j<a_i\}}+\sum_{j=1}^k \ind_{\{j<i,a_j=a_i\}}+1.
$$
The map $i\in\{1,\ldots,k\}\rightarrow R^\mathfrak{a}(i)\in\{1,\ldots,k\}$ is a bijection. Let $\rho^\mathfrak{a}$ be its inverse. Define $\uparrow\!\!\!\mathfrak{a}=\mathfrak{a}(\rho^\mathfrak{a})$, where $\uparrow\!\!\!\mathfrak{a}(j)$ can be roughly qualified as the $j$th value in the increasing order of $\{a_1,\ldots,a_k\}$. We check then that $\uparrow\!\!\!\mathfrak{a}$ is an non decreasing function on $\{1,\ldots,k\}$ taking the same values as $\mathfrak{a}$.

Let $1\leq j\leq k$, $i=\rho^\mathfrak{a}(j)$ and $b\in\mathbb{R}_+$. Let $\mathfrak{a}\nmid b=\{a_1\nmid b,\ldots,a_k\nmid b\}$. Suppose that $b\geq \uparrow\!\!\!\mathfrak{a}(j)=a_i$. Then, it can be checked that $a_h\nmid b<a_i\nmid b$ (resp. $a_h\nmid b=a_i\nmid b$) is equivalent to $a_h<a_i$ (resp. $a_h = a_i$). Therefore, $R^{\mathfrak{a}\nmid b}(i)=R^{\mathfrak{a}}(i) =j$ and $$
\uparrow\!\!\!(\mathfrak{a}\nmid b)(j)=\ a_i\nmid b = a_i = \ \uparrow\!\!\!\mathfrak{a}(j).
$$

\subsection{The enlargement of filtration with multiple random times}

Let $m>0$ be an integer and $\tau_1,\ldots,\tau_m$ be $m$ random times. For a $1\leq k\leq m$, consider the random function $\mathfrak{t}_k$ on $\{1,\ldots,k\}$ taking respectively the values $\{\tau_1,\ldots,\tau_k\}$. We define $\omega$ by $\omega$ the non decreasing function $\uparrow\!\!\!\mathfrak{t}_k$ as in the previous subsection.

\bl\label{re-ordering}
For any $1\leq j\leq k$, there exists a Borel function $\mathfrak{s}_{j}$ on $[0,\infty]^k$ such that $$
\uparrow\!\!\!\mathfrak{t}_k(j)=\mathfrak{s}_{j}(\tau_1,\ldots,\tau_k).
$$
If the $\tau_1,\ldots,\tau_k$ are stopping times with respect to some filtration, the random times $\uparrow\!\!\!\mathfrak{t}_k(1),\ldots,\uparrow\!\!\!\mathfrak{t}_k(k)$ also are stopping times with respect to the same filtration.
\el

\textbf{Proof.} This is a consequence of the following identity\!\!: for any $t\geq 0$,$$
\{\uparrow\!\!\!\mathfrak{t}_k(j)\leq t\}=\cup_{I\subset\{1,\ldots,k\}, \sharp I=j}\{\tau_h\leq t, \forall h\in I\}. \ok
$$

The random times $\uparrow\!\!\!\mathfrak{t}_k(1),\ldots,\uparrow\!\!\!\mathfrak{t}_k(k)$ form an increasing re-ordering of $\{\tau_1,\ldots,\tau_k\}$. We will denote them as $\sigma_{k,j}=\uparrow\!\!\!\mathfrak{t}_k(j), 1\leq j\leq k$.

Let $\mathbb{G}^0=\mathbb{F}$. For $1\leq k\leq m$, let $
\mathbb{G}^k=(\mathcal{G}^k_t)_{t\geq 0}$ where $$
\mathcal{G}^k_t=\mathcal{N}^k\vee(\cap_{s>t}(\mathcal{G}^{k-1}_s\vee\sigma(\tau_k\wedge s))),\ t\geq 0,
$$ 
and $$
\mathcal{N}^k
=\mathcal{N}^{\sigma(\tau_1)\vee\ldots\vee\sigma(\tau_k)\vee\mathcal{F}_\infty}.
$$
By induction, we can prove that $\mathbb{G}^k$ is the smallest right-continuous filtration containing $\mathbb{F}$ and $\mathcal{N}^k$, making the $\tau_1,\ldots,\tau_k$ stopping times.

\brem
Let us temporarily denote $\mathbb{G}^k$ by $\mathbb{G}^{(\tau_1,\ldots,\tau_k)}$ in reference to the dependence of $\mathbb{G}^k$ on the random times $(\tau_1,\ldots,\tau_k)$. We have the relation $$
\mathbb{G}^{(\tau_1,\ldots,\tau_k)}\supset\mathbb{G}^{(\sigma_{k,1},\ldots,\sigma_{k,k})},
$$
because the $\sigma_{k,j}, 1\leq j\leq k,$ are $\mathbb{G}^{(\tau_1,\ldots,\tau_k)}$-stopping times. In general, there is no equality between $\mathbb{G}^{(\sigma_{k,1},\ldots,\sigma_{k,k})}$ and $\mathbb{G}^{(\tau_1,\ldots,\tau_k)}$. For example, let $\{A,B,C\}$ be a partition of $\Omega$. If$$
\dcb
\tau_1=1\ind_A+2\ind_B+3\ind_C,\\
\tau_2=2\ind_A+3\ind_B+1\ind_C,\\
\tau_3=3\ind_A+1\ind_B+2\ind_C,\\
\dce
$$
we get $\sigma_{3,1}\equiv 1,\sigma_{3,2}\equiv 2,\sigma_{3,3}\equiv 3$.
\erem

\subsection{The optional splitting formulas}

\bd\label{df_1}
We say that the $\mathbb{G}^m$-optional splitting formula holds at times $\tau_1,\ldots,\tau_m$ with respect to $\mathbb{F}$, if, for any $\mathbb{G}^m$-optional process $Y$, there exist functions $Y^{(0)},Y^{(1)},\ldots, Y^{(m)}$ defined on $[0,\infty]^m\times(\mathbb{R}_+\times\Omega)$ being $\mathcal{B}[0,\infty]^m\otimes\mathcal{O}(\mathbb{F})$-measurable such that
$$
Y=\sum_{i=0}^mY^{(i)}(\tau_1\nmid\sigma_{m,i},\ldots,\tau_m\nmid\sigma_{m,i})\ind_{[\sigma_{m,i},\sigma_{m,i+1})}
$$
(an indistinguishable identity with respect to $\mathcal{N}^m$), where $\sigma_{m,0}\equiv 0$ and $\sigma_{m,m+1}\equiv \infty$ by definition.
\ed

Note that this definition is coherent with Definition \ref{df_splitting} when $m=1$.

\bl\label{parameter_formula}
Let $(E,\mathcal{E})$ be a measurable space. Let $1\leq k\leq m$. Suppose that the $\mathbb{G}^k$-optional splitting formula holds at times $\tau_1,\ldots,\tau_k$ with respect to $\mathbb{F}$. Then, for any $\mathcal{E}\otimes\mathcal{O}(\mathbb{G}^k)$-measurable function $Y(\theta,s,\omega)$, there exist functions $Y^{(0)},Y^{(1)},\ldots, Y^{(k)}$ defined on $E\times[0,\infty]^k\times(\mathbb{R}_+\times\Omega)$ being $\mathcal{E}\otimes\mathcal{B}[0,\infty]^k\otimes\mathcal{O}(\mathbb{F})$-measurable such that
$$
Y(\theta)=\sum_{i=0}^kY^{(i)}(\theta,\tau_1\nmid\sigma_{k,i},\ldots,\tau_k\nmid\sigma_{k,i})\ind_{[\sigma_{k,i},\sigma_{k,i+1})}.
$$
\el

\textbf{Proof} We only need to check the lemma upon the functions of the form $Y(\theta,s,\omega)=g(\theta)F_s(\omega), g\in\mathcal{E}, F\in\mathcal{O}(\mathbb{G}^k)$, and apply the monotone class theorem. \ok

\bethe\label{GmO}
Suppose $m>1$. Suppose that the $\mathbb{G}^{m-1}$-optional splitting formula holds at times $\tau_1,\ldots,\tau_{m-1}$ with respect to $\mathbb{F}$. Suppose the $\mathbb{G}^m$-optional splitting formula holds at time $\tau_m$ with respect to $\mathbb{G}^{m-1}$. Then, the $\mathbb{G}^m$-optional splitting formula holds at times $\tau_1,\ldots,\tau_{m-1},\tau_m$ with respect to $\mathbb{F}$.
\ethe

\textbf{Proof} Let $Y$ be a $\mathbb{G}^m$-optional process. By assumption, there exist $Y'$ and $Y''$ such that $Y'\in\mathcal{O}(\mathbb{G}^{m-1})$ and $Y''\in\mathcal{B}[0,\infty]\otimes\mathcal{O}(\mathbb{G}^{m-1})$ and $$
Y=Y'\ind_{[0,\tau_{m})}+Y''(\tau_{m})\ind_{[\tau_{m},\infty)}.
$$
The theorem will be proved, if we show that $Y'\ind_{[0,\tau_{m})}$ and $Y''(\tau_{m})\ind_{[\tau_{m},\infty)}$ satisfy the $\mathbb{G}^m$-optional splitting formula at $\tau_1,\ldots,\tau_m$ with respect to $\mathbb{F}$. To do this, we now rewrite the functions $Y'$ and $Y''(\tau_m)$ in terms of $\tau_h\nmid\sigma_{m,i}$.
 
According to the $\mathbb{G}^{m-1}$-optional splitting formula at times $\tau_1,\ldots,\tau_{m-1}$ with respect to $\mathbb{F}$, there exist functions $Y'^{(0)},Y'^{(1)},\ldots, Y'^{(m-1)}$ defined on $[0,\infty]^{m-1}\times(\mathbb{R}_+\times\Omega)$ being $\mathcal{B}[0,\infty]^{m-1}\otimes\mathcal{O}(\mathbb{F})$-measurable such that
$$
Y'=\sum_{i=0}^{m-1}Y'^{(i)}(\tau_1\nmid\sigma_{m-1,i},\ldots,\tau_{m-1}\nmid\sigma_{m-1,i})\ind_{[\sigma_{m-1,i},\sigma_{m-1,i+1})}.
$$
According to Lemma \ref{parameter_formula}, there exist functions $Y''^{(0)},Y''^{(1)},\ldots, Y''^{(m-1)}$ defined on $[0,\infty]\times[0,\infty]^{m-1}\times(\mathbb{R}_+\times\Omega)$ being $\mathcal{B}[0,\infty]\otimes\mathcal{B}[0,\infty]^{m-1}\otimes\mathcal{O}(\mathbb{F})$-measurable such that
$$
Y''(\theta)=\sum_{i=0}^{m-1}Y''^{(i)}(\theta,\tau_1\nmid\sigma_{m-1,i},\ldots,\tau_{m-1}\nmid\sigma_{m-1,i})\ind_{[\sigma_{m-1,i},\sigma_{m-1,i+1})}.
$$
The above expressions do not meet our requirements, because they employ $\sigma_{m-1,\cdot}$ instead of $\sigma_{m,\cdot}$. But a precise relationship exists which will make the transition from $\sigma_{m-1,\cdot}$ to $\sigma_{m,\cdot}$. If we denote $k=R^{\{\tau_1,\ldots,\tau_m\}}(m)$, we have $\sigma_{m-1,i}=\sigma_{m,i}$ for $i\leq k-1$ and $\sigma_{m-1,i}=\sigma_{m,i+1}$ for $k\leq i <m$. Moreover, if $k<m$, $\sigma_{m-1,k}=\sigma_{m,k+1}>\tau_m$. This relationship entails$$
\dcb
\ind_{[\sigma_{m-1,i},\ \sigma_{m-1,i+1})}\ind_{[0,\tau_{m})}
=\ind_{[\sigma_{m,i},\ \sigma_{m,i+1})}\ind_{[0,\sigma_{m,k})}
=\ind_{[\sigma_{m,i},\ \sigma_{m,i+1})},&&\mbox{ if $i\leq k-2$,}\\

\ind_{[\sigma_{m-1,k-1},\ \sigma_{m-1,k})}\ind_{[0,\tau_{m})}
=\ind_{[\sigma_{m,k-1},\ \sigma_{m,k+1})}\ind_{[0,\sigma_{m,k})}
=\ind_{[\sigma_{m,k-1},\ \sigma_{m,k})},&&\mbox{ if $i= k-1$,}\\

\ind_{[\sigma_{m-1,i},\ \sigma_{m-1,i+1})}\ind_{[0,\tau_{m})}=0,&&\mbox{ if $i\geq k$.}\\
\dce
$$
Notice that, if $i<k$ and $\sigma_{m,i}=\tau_m=\sigma_{m,k}$, necessarily $\sigma_{m,i+1}\leq \sigma_{m,k}$ so that $[\sigma_{m,i},\sigma_{m,i+1})=\emptyset$. 
Consequently,
$$
\dcb
\ind_{[\sigma_{m-1,i},\ \sigma_{m-1,i+1})}\ind_{[0,\tau_{m})}
=\ind_{[\sigma_{m,i},\ \sigma_{m,i+1})}
=\ind_{[\sigma_{m,i},\ \sigma_{m,i+1})}\ind_{\{\sigma_{m,i}<\tau_m\}},&&\mbox{ if $i\leq k-1$,}\\

\ind_{[\sigma_{m-1,i},\ \sigma_{m-1,i+1})}\ind_{[0,\tau_{m})}
=0=\ind_{[\sigma_{m,i},\ \sigma_{m,i+1})}\ind_{\{\sigma_{m,i}<\tau_m\}},&&\mbox{ if $i\geq k$.}\\
\dce
$$
With these identities we write
$$
\dcb
&&Y'\ind_{[0,\tau_{m})}\\
&=&\sum_{i=0}^{m-1}Y'^{(i)}(\tau_1\nmid\sigma_{m-1,i},\ldots,\tau_{m-1}\nmid\sigma_{m-1,i})\ind_{[\sigma_{m-1,i},\sigma_{m-1,i+1})}\ind_{[0,\tau_{m})}\\

&=&\sum_{i=0}^{k-1}Y'^{(i)}(\tau_1\nmid\sigma_{m-1,i},\ldots,\tau_{m-1}\nmid\sigma_{m-1,i})\ind_{[\sigma_{m-1,i},\sigma_{m-1,i+1})}\ind_{[0,\sigma_{m,k})}\\

&=&\sum_{i=0}^{k-1}Y'^{(i)}(\tau_1\nmid\sigma_{m,i},\ldots,\tau_{m-1}\nmid\sigma_{m,i})\ind_{[\sigma_{m,i},\sigma_{m,i+1})}\\

&=&\sum_{i=0}^{m-1}Y'^{(i)}(\tau_1\nmid\sigma_{m,i},\ldots,\tau_{m-1}\nmid\sigma_{m,i})\ind_{[\sigma_{m,i},\sigma_{m,i+1})}\ind_{\{\sigma_{m,i}<\tau_m\}}.
\dce
$$
In the last identity the condition $i\leq k-1$ is replaced by the condition $\sigma_{m,i}<\tau_m$. This is important because of the following relations $$
\dcb
\ind_{[\sigma_{m,i},\sigma_{m,i+1})}\ind_{\{\sigma_{m,i}<\tau_m\}}
&=&\ind_{[\sigma_{m,i},\sigma_{m,i+1})}\ind_{\{\sigma_{m,i}<\tau_m, \sigma_{m,i}<\infty\}}\\
&=&\ind_{[\sigma_{m,i},\sigma_{m,i+1})}\ind_{\{\tau_m\nmid\sigma_{m,i}=\infty, \sigma_{m,i}<\infty\}}\\
&=&\ind_{[\sigma_{m,i},\sigma_{m,i+1})}\ind_{\{\tau_m\nmid\sigma_{m,i}=\infty\}},
\dce
$$ 
which rewrites the expression with the term $\tau_m\nmid\sigma_{m,i}$. Substituting the last term into the preceding expression, we see that $Y'\ind_{[0,\tau_{m})}$ indeed satisfies the $\mathbb{G}^m$-optional splitting formula at $\tau_1,\ldots,\tau_m$ with respect to $\mathbb{F}$.

Next consider the interval $[\tau_{m},\infty)$. For any $j<k$, if $\sigma_{m,j}\geq \tau_m$, we necessarily have $\sigma_{m,j}=\sigma_{m,j+1}=\tau_m$, i.e., $[\sigma_{m,j},\sigma_{m,j+1})=\emptyset$. We compute
$$
\dcb
\ind_{[\sigma_{m-1,i},\ \sigma_{m-1,i+1})}\ind_{[\tau_{m},\infty)}
=\ind_{[\sigma_{m,i},\ \sigma_{m,i+1})}\ind_{[\sigma_{m,k},\infty)}
=0
=\ind_{[\sigma_{m,i+1},\ \sigma_{m,i+2})}\ind_{\{\tau_m\leq \sigma_{m,i+1}\}}
,&&\mbox{ if $i\leq k-2$,}\\

\ind_{[\sigma_{m-1,k-1},\ \sigma_{m-1,k})}\ind_{[\tau_{m},\infty)}
=\ind_{[\sigma_{m,k-1},\ \sigma_{m,k+1})}\ind_{[\sigma_{m,k},\infty)}
=\ind_{[\sigma_{m,k},\ \sigma_{m,k+1})}\ind_{\{\tau_m\leq \sigma_{m,k}\}},&&\mbox{ if $i= k-1$,}\\

\ind_{[\sigma_{m-1,i},\ \sigma_{m-1,i+1})}\ind_{[\tau_{m},\infty)}
=\ind_{[\sigma_{m,i+1},\ \sigma_{m,i+2})}\ind_{[\sigma_{m,k},\infty)}
=\ind_{[\sigma_{m,i+1},\ \sigma_{m,i+2})}\ind_{\{\tau_m\leq \sigma_{m,i+1}\}},&&\mbox{ if $i\geq k$,}\\
\dce
$$
and then
$$
\dcb
&&Y''(\tau_m)\ind_{[\tau_{m},\infty)}\\

&=&\sum_{i=0}^{m-1}Y''^{(i)}(\tau_m,\tau_1\nmid\sigma_{m-1,i},\ldots,\tau_{m-1}\nmid\sigma_{m-1,i})\ind_{[\sigma_{m-1,i},\sigma_{m-1,i+1})}\ind_{[\tau_{m},\infty)}\\

&=&\sum_{i=0}^{m-1}Y''^{(i)}(\tau_m,\tau_1\nmid\sigma_{m-1,i},\ldots,\tau_{m-1}\nmid\sigma_{m-1,i})\ind_{[\sigma_{m,i+1},\sigma_{m,i+2})}\ind_{\{\tau_m\leq \sigma_{m,i+1}\}}\\

&=&\sum_{j=1}^{m}Y''^{(j-1)}(\tau_m\nmid\sigma_{m,j},\tau_1\nmid\sigma_{m-1,j-1},\ldots,\tau_{m-1}\nmid\sigma_{m-1,j-1})\ind_{[\sigma_{m,j},\sigma_{m,j+1})}\ind_{\{\tau_{m}\leq \sigma_{m,j}\}}.\\

\dce
$$
Notice that we have not directly substituted $\sigma_{m-1,\cdot}$ with $\sigma_{m,\cdot}$ in the expression 
$$
Y''^{(j-1)}(\tau_m\nmid\sigma_{m,j},\tau_1\nmid\sigma_{m-1,j-1},\ldots,\tau_{m-1}\nmid\sigma_{m-1,j-1})\ind_{[\sigma_{m,j},\sigma_{m,j+1})}\ind_{\{\tau_{m}\leq \sigma_{m,j}\}}.
$$
This is because, according to whether or not $j\leq k$, the substitutes are different. Since $k$ is random, we have to be careful about the measurability issue of such a substitution. We begin with $$
\ind_{[\sigma_{m,j},\sigma_{m,j+1})}\ind_{\{\tau_m\leq \sigma_{m,j}\}}
=\ind_{[\sigma_{m,j},\sigma_{m,j+1})}\ind_{\{\tau_m\nmid\sigma_{m,j}<\infty\}}.
$$ 
We continue with $\tau_a\nmid\sigma_{m-1,j-1}$ for $1\leq a<m$ and $1\leq j\leq m$. Since $\sigma_{m,j}\geq \sigma_{m-1,j-1}$, we have$$
\tau_a\nmid\sigma_{m-1,j-1}=(\tau_a\nmid\sigma_{m,j})\nmid\sigma_{m-1,j-1},
$$
and, according to subsection \ref{ordering}, $$
\sigma_{m-1,j-1}=\uparrow\!\!\!(\mathfrak{t}_{m-1}\nmid \sigma_{m,j})(j-1).
$$
By Lemma \ref{re-ordering}, there exists a function $\mathfrak{s}_{m-1,j-1}$ such that $$
\sigma_{m-1,j-1}=\uparrow\!\!\!(\mathfrak{t}_{m-1}\nmid \sigma_{m,j})(j-1)
=\mathfrak{s}_{m-1,j-1}(\tau_1\nmid \sigma_{m,j},\ldots,\tau_{m-1}\nmid \sigma_{m,j}).
$$
Consequently, $\tau_a\nmid\sigma_{m-1,j-1}$ also is a Borel function of $(\tau_1\nmid \sigma_{m,j},\ldots,\tau_{m-1}\nmid \sigma_{m,j})$.
From these facts, we conclude that there exist Borel functions $Z''^{(j-1)}$ on $[0,\infty]^{m}\times(\mathbb{R}_+\times\Omega)$ being $\mathcal{B}[0,\infty]^{m}\otimes\mathcal{O}(\mathbb{F})$-measurable such that 
$$
\dcb
&&Y''^{(j-1)}(\tau_m\nmid\sigma_{m,j},\tau_1\nmid\sigma_{m-1,j-1},\ldots,\tau_{m-1}\nmid\sigma_{m-1,j-1})\ind_{[\sigma_{m,j},\sigma_{m,j+1})}\ind_{\{\tau_{m}\leq \sigma_{m,j}\}}\\
&=&
Z''^{(j-1)}(\tau_1\nmid\sigma_{m,j},\ldots,\tau_{m-1}\nmid\sigma_{m,j}, \tau_m\nmid\sigma_{m,j})\ind_{[\sigma_{m,j},\sigma_{m,j+1})}.
\dce
$$
Substituting $Z''^{(j-1)}$ into the expression of $Y''(\tau_m)\ind_{[\tau_{m},\infty)}$, we finally prove that $Y''(\tau_m)\ind_{[\tau_{m},\infty)}$ satisfies the $\mathbb{G}^m$-optional splitting formula at times $\tau_1,\ldots,\tau_m$ with respect to $\mathbb{F}$. \ok

\subsection{Density hypothesis}

\bd\label{density}
We say that $(\tau_1,\ldots,\tau_m)$ satisfies the (conditional) density hypothesis with respect to $\mathcal{F}_\infty$, if there exists a Borel probability measure $\mu$ on $[0,\infty]$ and a non negative function $\gamma$ on $[0,\infty]^m\times\Omega$ being $\mathcal{B}[0,\infty]^m\otimes\mathcal{F}_\infty$ measurable such that$$
\mathbb{Q}[(\tau_1,\ldots,\tau_m)\in\mathtt{A}\ |\mathcal{F}_\infty]
=\int_\mathtt{A}\gamma(t_1,\ldots,t_m)\mu^{\otimes m}(\mathsf{d}t_1,\ldots,\mathsf{d}t_m)
$$
for any $\mathtt{A}\in\mathcal{B}[0,\infty]^m$.
\ed

We have the following results.

\bl\label{density-density}
Suppose that $(\tau_1,\ldots,\tau_m)$ satisfies the density hypothesis with respect to $\mathcal{F}_\infty$. Then, $\tau_m$ satisfies the density hypothesis with respect to $\mathcal{G}^{m-1}_\infty$. For any $1\leq k< m$, $(\tau_1,\ldots,\tau_k)$ satisfies the density hypothesis with respect to $\mathcal{F}_\infty$.
\el

The proof of the lemma is straightforward. 

\bl\label{one-tau-density}
Consider the case of $m=1$. If the density hypothesis holds for $\tau_1$ with respect to $\mathcal{F}_\infty$, the global $\mathbb{G}^1$-optional splitting formula holds at $\tau_1$ with respect to $\mathbb{F}$.
\el

\textbf{Proof} 
The proof is based on Lemma \ref{LSST}. Let $h(u,\omega)$ be a bounded function defined on $[0,\infty]\times\Omega$, $\mathcal{B}[0,\infty]\otimes\mathcal{F}_\infty$ measurable. A direct computation with the density hypothesis yields$$
\dcb
&&\mathbb{E}[h(\tau)|\mathcal{F}_{t}\vee\sigma(\tau\nmid t)]\\
&=&\ind_{\{t<\tau\}}\frac{\mathbb{E}[h(\tau)\ind_{\{t<\tau\}}|\mathcal{F}_{t}]}{\mathbb{E}[\ind_{\{t<\tau\}}|\mathcal{F}_{t}]}
+
\ind_{\{\tau\leq t\}}\left.\frac{\mathbb{E}[h(u)\gamma(u)|\mathcal{F}_{t}]}{\mathbb{E}[\gamma(u)|\mathcal{F}_{t}]}\ind_{\{\mathbb{E}[\gamma(u)|\mathcal{F}_{t}]>0\}}\right|_{u=\tau},
\dce
$$
where $\mathbb{E}[h(u)\gamma(u)|\mathcal{F}_{t}]$ denotes the value at $t$ of the parametered $\mathbb{F}$ optional projection of the parametered random variable $h(u,\omega)\gamma(u,\omega)$, introduced in \cite[Proposition 3]{SY} (similar interpretation for the notation $\mathbb{E}[\gamma(u)|\mathcal{F}_{t}]$). Notice that, by \cite[Chapitre VI n$^\circ$48]{DM}, these parametered $\mathbb{F}$ optional projections have right continuous path. Taking the right limit in the above formula we see that the martingale $\mathbb{E}[h(\tau)|\mathcal{G}_{t}], t\in\mathbb{R}_+$, satisfies the optional splitting formula.

Notice that $\mathcal{G}_{\infty-}$ is generated by $\mathcal{F}_{\infty-}\vee\sigma(\tau)$ together with negligible sets. Hence, any bounded $\mathbb{G}$ martingale $X$ is indistinguishable to a martingale of the type considered in the previous paragraph. The lemma is proved because of Lemma \ref{LSST}. \ok

Now look at Lemma \ref{one-tau-density}, Theorem \ref{GmO}, and Lemma \ref{density-density}. They constitute a perfect mathematical induction pattern. We obtain the following result\!\!:

\bethe\label{densityH-splitting}
If times $(\tau_1,\ldots,\tau_m)$ satisfy the density hypothesis with respect to $\mathcal{F}_\infty$, then, $\mathbb{G}^m$-optional splitting formula holds at times $\tau_1,\ldots,\tau_m$ with respect to $\mathbb{F}$. 
\ethe

\

\section{Splitting formula at random times with marks}\label{multi-time-marked}

The results in the previous section can be extended to the case of random times with marks. The proofs follow the same idea with some notational complications.

\subsection{Filtration $\mathbb{G}^{*m}$ and optional splitting formula}

Let $(E,\mathcal{E})$ be a separable complete metric space with its Borel $\sigma$-algebra. Let $\vartriangle\in E$ and $E^\circ=E\setminus \{\vartriangle\}$. Let $(\xi_1,\ldots,\xi_m)$ be $m$ random variables taking values in $E^\circ$. Define, for $1\leq i\leq m, t\geq 0$,$$
H_{i}(t)=
\left\{
\dcb
\vartriangle&&\mbox{if $t <\tau_i$,}\\
\\
\xi_i&&\mbox{if $\tau_i\leq t$.}
\dce
\right.\ \ \  
$$
Let $\mathcal{H}^{\{1,\ldots,m\}}_t=\sigma(H_i(s)\!\!: 1\leq i\leq m,\ 0\leq s\leq t)$ and $$
\mathcal{G}^{*m}_t=\mathcal{N}^{*m}\vee\cap_{s>t}(\mathcal{F}_s\vee\mathcal{H}^{\{1,\ldots,m\}}_s),
$$
where $\mathcal{N}^{*m}$ denotes $\mathcal{N}^{\mathcal{H}^{\{1,\ldots,m\}}_\infty\vee\mathcal{F}_\infty}$. Let $\mathbb{G}^{*m}$ be the filtration of $\mathcal{G}^{*m}_t, t\geq 0$. 

Let $\mathfrak{D}(E)$ be the space of all càdlàg functions taking values in $E$ equipped with the Skorokhod topology (cf. \cite{billingsley}) and its Borel $\sigma$-algebra $\mathcal{D}$.

\bd\label{df-marked-times}
We say that the $\mathbb{G}^{*m}$-optional splitting formula holds at times $\tau_1,\ldots,\tau_m$ with respect to $\mathbb{F}$, if, for any $\mathbb{G}^{*m}$-optional process $Y$, there exist functions $Y^{(0)},Y^{(1)},\ldots, Y^{(m)}$ defined on $\mathfrak{D}(E)^m\times(\mathbb{R}_+\times\Omega)$ being $\mathcal{D}^{m}\otimes\mathcal{O}(\mathbb{F})$-measurable such that
$$
Y=\sum_{i=0}^mY^{(i)}(H_1^{\sigma_{m,i}},\ldots,H_m^{\sigma_{m,i}})\ind_{[\sigma_{m,i},\sigma_{m,i+1})},
$$
where $H_i^{\sigma_{m,i}}$ denotes the process $H_i$ stopped at $\sigma_{m,i}$.
\ed

Note that, for $1\leq k\leq m, 0\leq i\leq m, 0\leq u<\infty$, $$
\dcb
H_k^{\sigma_{m,i}}(u)
&=&\ind_{\{\sigma_{m,i}<\tau_k\mbox{ or }u<\tau_k\}}\Delta 
+\ind_{\{\sigma_{m,i}\geq \tau_k, u\geq \tau_k\}}\xi_k\\

&=&\ind_{\{\tau_k\nmid\sigma_{m,i}> u\}}\Delta 
+\ind_{\{\tau_k\nmid\sigma_{m,i}\leq u \}}\xi_k.\\
\dce
$$
So, if $\xi_k$ are constant random variables, the above Definition \ref{df-marked-times} coincides with Definition \ref{df_1}.

\bethe\label{GmO-star}
Suppose $m>1$. Supose that $\mathbb{G}^{*m-1}$-optional splitting formula holds at times $\tau_1,\ldots,\tau_{m-1}$ with respect to $\mathbb{F}$. Suppose $\mathbb{G}^{*m}$-optional splitting formula holds at time $\tau_m$ with respect to $\mathbb{G}^{*m-1}$. Then, $\mathbb{G}^{*m}$-optional splitting formula holds at times $\tau_1,\ldots,\tau_m$ with respect to $\mathbb{F}$.
\ethe

\textbf{Proof} Let $Y$ be a $\mathbb{G}^{*m}$-optional process. By assumption, there exist $Y',Y''$ such that $Y'\in\mathcal{O}(\mathbb{G}^{*m-1})$ and $Y''\in\mathcal{D}\otimes\mathcal{O}(\mathbb{G}^{*m-1})$ and $$
Y=Y'\ind_{[0,\tau_{m})}+Y''(H_m)\ind_{[\tau_{m},\infty)}.
$$
Now, according to the $\mathbb{G}^{*m-1}$-optional splitting formula at times $\tau_1,\ldots,\tau_{m-1}$ with respect to $\mathbb{F}$, there exist functions $Y'^{(0)},Y'^{(1)},\ldots, Y'^{(m-1)}$ defined on $\mathfrak{D}(E)^{m-1}\times(\mathbb{R}_+\times\Omega)$ being $\mathcal{D}^{m-1}\otimes\mathcal{O}(\mathbb{F})$-measurable such that
$$
Y'=\sum_{i=0}^{m-1}Y'^{(i)}(H_1^{\sigma_{m-1,i}},\ldots,H_{m-1}^{\sigma_{m-1,i}})\ind_{[\sigma_{m-1,i},\sigma_{m-1,i+1})}.
$$
Also, there exist functions $Y''^{(0)},Y''^{(1)},\ldots, Y''^{(m-1)}$ defined on $\mathfrak{D}(E)\times\mathfrak{D}(E)^{m-1}\times(\mathbb{R}_+\times\Omega)$ being $\mathcal{D}\otimes\mathcal{D}^{m-1}\otimes\mathcal{O}(\mathbb{F})$-measurable such that
$$
Y''(H_m)=\sum_{i=0}^{m-1}Y''^{(i)}(H_m,H_1^{\sigma_{m-1,i}},\ldots,H_{m-1}^{\sigma_{m-1,i}})\ind_{[\sigma_{m-1,i},\sigma_{m-1,i+1})}.
$$
We have (cf. the proof of Theorem \ref{GmO})
$$
\dcb
Y'\ind_{[0,\tau_{m})}

&=&\sum_{i=0}^{m-1}Y'^{(i)}(H_1^{\sigma_{m,i}},\ldots,H_{m-1}^{\sigma_{m,i}})\ind_{[\sigma_{m,i},\sigma_{m,i+1})}\ind_{\{\sigma_{m,i}<\tau_m\}}.
\dce
$$
Since $\ind_{\{\sigma_{m,i}<\tau_m\}}
=\ind_{\{H_m(\sigma_{m,i})=\vartriangle\}}
=\ind_{\{H_m^{\sigma_{m,i}}(\infty)=\vartriangle, \sigma_{m,i}<\infty\}}$, 
the above expression is simply a $\mathbb{G}^{*m}$-optional splitting formula for $Y'\ind_{[0,\tau_{m})}$ at $\tau_1,\ldots,\tau_m$ with respect to $\mathbb{F}$.

Next, we write (cf. the proof of Theorem \ref{GmO})\!\!:
$$
\dcb
&&Y''(H_m)\ind_{[\tau_{m},\infty)}\\

&=&\sum_{j=1}^{m}Y''^{(j-1)}(H_m,H_1^{\sigma_{m-1,j-1}},\ldots,H_{m-1}^{\sigma_{m-1,j-1}})\ind_{[\sigma_{m,j},\sigma_{m,j+1})}\ind_{\{\tau_m\leq \sigma_{m,j}\}}\\

&=&\sum_{j=1}^{m}Y''^{(j-1)}(H_m^{\sigma_{m,j}},H_1^{\sigma_{m-1,j-1}},\ldots,H_{m-1}^{\sigma_{m-1,j-1}})\ind_{[\sigma_{m,j},\sigma_{m,j+1})}\ind_{\{H_m(\sigma_{m,j})\neq \vartriangle\}}.\\
\dce
$$
Recall that $$
\sigma_{m-1,j-1}=\uparrow\!\!\!(\mathfrak{t}_{m-1}\nmid \sigma_{m,j})(j-1)
=\mathfrak{s}_{m-1,j-1}(\tau_1\nmid \sigma_{m,j},\ldots,\tau_{m-1}\nmid \sigma_{m,j}), \ 1\leq i\leq m-1,
$$
and also, for $1\leq a\leq m-1$ and $u\geq 0$,
$$
\dcb
H_a^{\sigma_{m-1,j-1}}(u)
&=&(H_a^{\sigma_{m,j}})^{\sigma_{m-1,j-1}}(u)
=H_a^{\sigma_{m,j}}(u\wedge \sigma_{m-1,j-1}).\\
\dce
$$
Consider the measurability of the above object. First of all,$$
\forall 0\leq t<\infty,\ \{\tau_a\nmid\sigma_{m,j}\leq t\}
=\{\tau_a\leq t, \tau_a\leq\sigma_{m,j}\}
=\{H_a^{\sigma_{m,j}}(t)\in E^\circ\}
\in\sigma(H_a^{\sigma_{m,j}}).
$$
This implies that $$
\sigma_{m-1,j-1}\in\sigma(H_1^{\sigma_{m,j}},\ldots,H_{m-1}^{\sigma_{m,j}}).
$$
Note that the map $(t,\omega)\longrightarrow H_a^{\sigma_{m,j}}(t)$ is $\mathcal{B}[0,\infty]\otimes\sigma(H_a^{\sigma_{m,j}})$ measurable. Composing this map with that one $\omega\longrightarrow (u\wedge \sigma_{m-1,j-1}(\omega),\omega)$, we obtain that $$
H_a^{\sigma_{m-1,i}}(u)=H_a^{\sigma_{m,j}}(u\wedge \sigma_{m-1,j-1})
\in\sigma(H_1^{\sigma_{m,j}},\ldots,H_{m-1}^{\sigma_{m,j}}).
$$
Consequently, there exist Borel functions $Z''^{(i)}$ on $\mathfrak{D}(E)^{m}\times(\mathbb{R}_+\times\Omega)$ being $\mathcal{D}^{m}\otimes\mathcal{O}(\mathbb{F})$-measurable such that 
$$
\dcb
&&Y''(H_m)\ind_{[\tau_{m},\infty)}\\

&=&\sum_{j=1}^{m}Y''^{(j-1)}(H_m^{\sigma_{m,j}},H_1^{\sigma_{m-1,j-1}},\ldots,H_{m-1}^{\sigma_{m-1,j-1}})\ind_{[\sigma_{m,j},\sigma_{m,j+1})}\ind_{\{H_m(\sigma_{m,j})\neq \vartriangle\}}\\

&=&\sum_{j=1}^{m}Z''^{(j)}(H_1^{\sigma_{m,j}},\ldots,H_{m-1}^{\sigma_{m,j}},H_m^{\sigma_{m,j}})\ind_{[\sigma_{m,j},\sigma_{m,j+1})}.\\
\dce
$$
This is a $\mathbb{G}^{*m}$-optional splitting formula for $Y''(\tau_m)\ind_{[\tau_{m},\infty)}$ at times $\tau_1,\ldots,\tau_m$ with respect to $\mathbb{F}$. \ok

\subsection{Density hypothesis}

\bd\label{density-marked}
We say that $((\xi_1,\tau_1),\ldots,(\xi_m,\tau_m))$ satisfies the (conditional) density hypothesis with respect to $\mathcal{F}_\infty$, if there exists a Borel probability measure $\nu$ on $E\times[0,\infty]$ and a non negative function $\gamma^*$ on $(E\times[0,\infty])^m\times\Omega$ being $(\mathcal{E}\otimes\mathcal{B}[0,\infty])^m\otimes\mathcal{F}_\infty$ measurable such that$$
\mathbb{Q}[((\xi_1,\tau_1),\ldots,(\xi_m,\tau_m))\in\mathtt{A}\ |\mathcal{F}_\infty]
=\int_\mathtt{A}\gamma^*((x_1,t_1),\ldots,(x_m,t_m))\nu^{\otimes m}(\mathsf{d}(x_1,t_1),\ldots,\mathsf{d}(x_m,t_m))
$$
for any $\mathtt{A}\in(\mathcal{E}\times\mathcal{B}[0,\infty])^m$.
\ed

Notice that the probability measure $\nu$ in Definition \ref{density-marked} necessarily has support $E^\circ\times[0,\infty]$.

\bl\label{one-tau-density-marked}
If $m=1$, if $(\xi_1,\tau_1)$ satisfies the density hypothesis with respect to $\mathcal{F}_\infty$, then the $\mathbb{G}^{*1}$-optional splitting formula holds at $\tau_1$ with respect to $\mathbb{F}$.
\el

\textbf{Proof} 
Let $t\in\mathbb{R}_+$. Let $h(x,u,\omega)$ (respectively $f(x,u,\omega)$) be a non negative function defined on $E\times[0,\infty]\times\Omega$, $\mathcal{E}\otimes\mathcal{B}[0,\infty]\otimes\mathcal{F}_\infty$ measurable (respectively $\mathcal{E}\otimes\mathcal{B}[0,\infty]\otimes\mathcal{F}_t$ measurable). For $u\in[0,\infty], x\in E$, denote by $x^{(u)}$ the point $x$ if $u<\infty$ or the point $\Delta$ if $u=\infty$. Denote by $\mathbb{E}[h(x,u)\gamma^*(x,u)|\mathcal{F}_{t}]$ the value at $t$ of the parametered $\mathbb{F}$ optional projection of the random variable $h(x,u,\omega)\gamma^*(x,u,\omega)$ with parameter $(x,u)$, introduced in \cite[Proposition 3]{SY} (similar interpretation for the notation $\mathbb{E}[\gamma^*(x,u)|\mathcal{F}_{t}]$). We compute$$
\dcb
&&\mathbb{E}[h(\xi_1,\tau_1)\ind_{\{t<\tau_1\}}f(\xi_1^{(\tau_1\nmid t)},\tau_1\nmid t)]\\
&=&\mathbb{E}[h(\xi_1,\tau_1)\ind_{\{t<\tau_1\}}f(\Delta,\infty)]\\

&=&\mathbb{E}[\ \mathbb{E}[h(\xi_1,\tau_1)\ind_{\{t<\tau_1\}}|\mathcal{F}_t]\ f(\Delta,\infty)]\\

&=&\mathbb{E}[\frac{\mathbb{E}[h(\xi_1,\tau_1)\ind_{\{t<\tau_1\}}|\mathcal{F}_t]}{\mathbb{E}[\ind_{\{t<\tau_1\}}|\mathcal{F}_t]}\ind_{\{\mathbb{E}[\ind_{\{t<\tau_1\}}|\mathcal{F}_t]>0\}} \ind_{\{t<\tau_1\}} f(\xi_1^{(\tau_1\nmid t)},\tau_1\nmid t)],
\dce
$$
and
$$
\dcb
&&\mathbb{E}[h(\xi_1,\tau_1)\ind_{\{\tau_1\leq t\}}f(\xi_1^{(\tau_1\nmid t)},\tau_1\nmid t)]\\

&=&\mathbb{E}[h(\xi_1,\tau_1)\ind_{\{\tau_1\leq t\}}f(\xi_1,\tau_1)]\\

&=&\mathbb{E}[\int h(x,u)\ind_{\{u\leq t\}}f(x,u)\gamma^*(x,u)\nu(dx,du)]\\

&=&\int\mathbb{E}[ h(x,u)\ind_{\{u\leq t\}}f(x,u)\gamma^*(x,u)]\nu(dx,du)\\

&=&\int \mathbb{E}[\ \mathbb{E}[h(x,u)\gamma^*(x,u)|\mathcal{F}_t]\ind_{\{u\leq t\}} f(x,u)]\nu(dx,du)\\

&=&\int \mathbb{E}[\frac{\mathbb{E}[h(x,u)\gamma^*(x,u)|\mathcal{F}_t]}{\mathbb{E}[\gamma^*(x,u)|\mathcal{F}_t]}\ind_{\{\mathbb{E}[\gamma^*(x,u)|\mathcal{F}_t]>0\}}\ind_{\{u\leq t\}} f(x,u)\gamma^*(x,u)]\nu(dx,du)\\

&=& \mathbb{E}[\int\frac{\mathbb{E}[h(x,u)\gamma^*(x,u)|\mathcal{F}_t]}{\mathbb{E}[\gamma^*(x,u)|\mathcal{F}_t]}\ind_{\{\mathbb{E}[\gamma^*(x,u)|\mathcal{F}_t]>0\}}\ind_{\{u\leq t\}} f(x,u)\gamma^*(x,u)\nu(dx,du)]\\

&=& \mathbb{E}[\left(\frac{\mathbb{E}[h(x,u)\gamma^*(x,u)|\mathcal{F}_t]}{\mathbb{E}[\gamma^*(x,u)|\mathcal{F}_t]}\ind_{\{\mathbb{E}[\gamma^*(x,u)|\mathcal{F}_t]>0\}}\right)_{x=\xi_1,u=\tau_i} \ind_{\{\tau_1\leq t\}} f(\xi_1^{(\tau_1\nmid t)},\tau_1\nmid t)].
\dce
$$
This computation shows$$
\dcb
&&\mathbb{E}[h(\xi_1,\tau_1)|\mathcal{F}_t\vee\sigma(\xi_1^{(\tau_1\nmid t)},\tau_1\nmid t)]\\
&=&
\ind_{\{t<\tau_1\}}\frac{\mathbb{E}[h(\xi_1,\tau_1)\ind_{\{t<\tau_1\}}|\mathcal{F}_t]}{\mathbb{E}[\ind_{\{t<\tau_1\}}|\mathcal{F}_t]}\ind_{\{\mathbb{E}[\ind_{\{t<\tau_1\}}|\mathcal{F}_t]>0\}}
+
\ind_{\{\tau_1\leq t\}}\left(\frac{\mathbb{E}[h(x,u)\gamma^*(x,u)|\mathcal{F}_t]}{\mathbb{E}[\gamma^*(x,u)|\mathcal{F}_t]}\ind_{\{\mathbb{E}[\gamma^*(x,u)|\mathcal{F}_t]>0\}}\right)_{x=\xi_1,u=\tau_i}.
\dce
$$
Notice that, by \cite[Chapitre VI n$^\circ$48]{DM}, the parametered $\mathbb{F}$ optional projections in this formula have right continuous path. Taking the right limit in the above formula, we see that the martingale $\mathbb{E}[h(\xi_1,\tau_1)|\mathcal{G}_{t}], t\in\mathbb{R}_+$, satisfies the optional splitting formula.

Notice that $\mathcal{G}_{\infty-}$ is generated by $\mathcal{F}_{\infty-}\vee\sigma(\xi_1,\tau_1)$ together with negligible sets. Hence, any bounded $\mathbb{G}$ martingale $X$ is indistinguishable to a martingale of the type considered in the previous paragraph. Now to complete the proof of the lemma, we only need to repeat the argument in \cite[Chapitre XX, n$^\circ$22]{DM2}, as we did at the end of the proof of Lemma \ref{LSST}. \ok

Reproducing the argument in the proof of Theorem \ref{densityH-splitting}, we also obtain\!\!:

\bethe\label{densityH-splitting-marked}
If the marked times $((\xi_1,\tau_1),\ldots,(\xi_m,\tau_m))$ satisfy the density hypothesis with respect to $\mathcal{F}_\infty$, then, $\mathbb{G}^{*m}$-optional splitting formula holds at times $\tau_1,\ldots,\tau_m$ with respect to $\mathbb{F}$. 
\ethe

\brem
As a consequence of the above theorem, the optional splitting formula is valid in the papers \cite{CJZ, ElKJJ, jiao, KL, Pham} due to the density hypothesis.
\erem

\

\textbf{Acknowledgment} I am thankful to Lim T. and Rutkowski M. for the discussions I had with them about the optional splitting formula. This research benefited from the support of the "Chair Markets in Transition", under the aegis of Louis Bachelier laboratory, a joint initiative of Ecole polytechnique, Université d'Evry Val d'Essonne and Fédération Bancaire Française.

\

\end{document}